\documentclass[aap,seceqn,MSNbibl,citesort,dvips]{arximspdf}

\doi{10.1214/10-AAP717}
\volume{21}
\issue{4}
\pubyear{2011}
\firstpage{1282}
\lastpage{1321}

\makeatletter

\newcommand{\varliminf}{\mathop{\underline{\lim}}}
\newcommand{\varlimsup}{\mathop{\overline{\lim}}}

\newcommand{\R}{\mathbb R}
\newcommand{\Prob}{\mathbb{P}}
\newcommand{\Hh}{\mathcal{H}}
\newcommand{\esp}{\mathbb{E}}
\newcommand{\s}{\sigma_{X_t}}
\newcommand{\dA}{\partial A}
\newcommand{\nB}{|B|}
\newcommand{\nA}{|A|}
\newcommand{\DX}{D^{j_1,\ldots,j_k}_{r_1,\ldots,r_k} X^i_t}
\newcommand{\Bun}{B_R (y_0)}
\newcommand{\Bdeux}{B_{3 R} (y_0) }
\newcommand{\pcap}{\widehat{p}_{t,y_0}}
\newcommand{\pt}{p_{t,y_0}}
\newcommand{\nut}{\nu_{t,\delta}}
\newcommand{\taut}{\tau_{t, \delta}}
\newcommand{\Xbar}{\overline{X}}
\newcommand{\gammaHat}{\widehat{\gamma}}

\newproclaim{defin}{Definition}[section]
\newtheorem{theor}{Theorem}[section]
\newtheorem{lem}{Lemma}[section]
\newtheorem{prop}{Proposition}[section]
\newproclaim{rem}{Remark}[section]
\newtheorem{cor}{Corollary}

\def\bsuffix #1{#1}

\makeatother

\begin{document}
\begin{frontmatter}

\title{Smoothness and asymptotic estimates of densities
for SDEs with locally smooth coefficients and applications to square
root-type diffusions}
\runtitle{Smoothness of densities and asymptotic estimates}

\begin{aug}
\author[A]{\fnms{Stefano} \snm{De Marco}\corref{}\thanksref{t1}\ead[label=e1]{s.demarco@sns.it}}
\runauthor{S. De Marco}
\affiliation{Universit\'e Paris-Est and Scuola
Normale Superiore
di Pisa}
\address[A]{Laboratoire d'Analyse et Math\'ematiques\\
\quad Appliqu\'ees UMR 8050\\
Universit\'e Paris-Est\\
5 bd Descartes, Champs-sur-Marne\\
77454 Marne-la-Vall\'ee Cedex 2\\
France\\
and\\
Scuola Normale Superiore di Pisa\\
Piazza dei Cavalieri 7\\
56126 Pisa\\
Italy\\
\printead{e1}} 
\end{aug}

\thankstext{t1}{Supported by the European Science Foundation
(ESF) within the framework of the ESF activity entitled ``Advanced
Mathematical Methods for
Finance'' for the first semester of 2009.}

\received{\smonth{4} \syear{2009}}
\revised{\smonth{5} \syear{2010}}

%
\begin{abstract}
We study smoothness of densities for the solutions of SDEs whose
coefficients are smooth and nondegenerate only on an open domain~$D$.
We prove that a smooth density exists on $D$ and give upper bounds~for
this density. Under some additional conditions (mainly dealing with the
growth of the coefficients and their derivatives), we formulate upper
bounds that are suitable to obtain asymptotic estimates of the density
for large values of the state variable (``tail'' estimates). These
results specify and extend some results by Kusuoka and Stroock
[\textit{J.~Fac. Sci. Univ. Tokyo Sect. IA Math.}
\textbf{32} (1985) 1--76], but our approach is substantially different
and based on a technique to estimate the Fourier transform inspired
from Fournier [\textit{Electron. J.~Probab.} \textbf{13} (2008) 135--156]
and Bally [Integration by parts formula for locally smooth laws and applications
to equations with jumps I (2007) The Royal Swedish Academy of
Sciences]. This
study is motivated by existing models for financial securities which
rely on SDEs with non-Lipschitz coefficients. Indeed, we apply our
results to a square root-type diffusion (CIR or CEV) with coefficients
depending on the state variable, that is, a situation where standard
techniques for density estimation based on Malliavin calculus do not
apply. We establish the existence of a smooth density, for which we
give exponential estimates and study the behavior at the origin (the
singular point).
\end{abstract}

%

\end{frontmatter}

\section{Introduction}

It is well known that
Malliavin calculus is a tool which allows, among other, to prove that
the law of a diffusion process admits a~smooth density.
More precisely, if one assumes that the coefficients of an SDE are
bounded $C^{\infty}$ functions
with bounded derivatives of any order and that, on the other hand, the
Hormand\"{e}r condition holds,
then the solution of the equation is a smooth functional in Malliavin's sense, and it
is nondegenerate at any fixed
positive time.
Then the general criterion given by Malliavin \cite{Mall1} 
allows one to say that the law of such a random variable is absolutely
continuous with respect to the
Lebesgue measure, and its density is a smooth function (see \cite
{Nual} for a general presentation of this topic).

The aim of this paper is to relax the aforementioned conditions on the
coefficients: roughly
speaking, we assume that the coefficients are smooth only on an open
domain $D$ and have bounded
partial derivatives therein.
Moreover, we assume that the nondegeneracy condition on the diffusion
coefficient
holds true on $D$ only.
Under these assumptions, we prove that the law of a strong solution to
the equation admits a smooth density on $D$ (Theorem \ref{t:localDens1}).
Furthermore, when $D$ is the complementary of a compact ball and the
coefficients satisfy some
additional assumptions on $D$ (mainly dealing with their growth and the
one of derivatives), we give
upper bounds for the density for large values of the state variable
(Theorem~\ref{t:localDensSubLin1}).
We will occasionally refer to these aymptotic estimates of the density
as ``tail estimates'' or estimates on the density's ``tails.''

Local results have already been obtained by Kusuoka and Stroock
in \cite{KusSt}, Section 4.
Here the authors work under local regularity and nondegeneracy
hypotheses too,
but the bounds they provide on the density are mostly significant
on the diagonal (i.e., close to starting point) and in the small time limit,
while they are not appropriate for tail estimates.
Moreover, the constants appearing in the estimates are not explicit
(cf. (4.7)--(4.9) in Theo\-rem~4.5 and the
corresponding estimates in Corollary 4.10, \cite{KusSt}).
In the present paper, we provide upper bounds that are suitable for tail
estimates, and we find out the
explicit dependence of the bounding constants with respect to the
coefficients of the SDE and their
derivatives.
Our bounds turn out to be applicable to the case of diffusions with
tails stronger than gaussian.
This is the case for square-root diffusions, which
are our major example of interest
(see Section \ref{s:localCIR}).
Also, our approach is substantially different from the one in \cite{KusSt}.
In particular, we rely on a Fourier transform argument, employing a
technique to estimate the
Fourier transform of the process inspired from the work of Fournier in
\cite{Fournier}
and of Bally in \cite{BalLocSmoothLaws} and relying on
specifically-designed Malliavin calculus techniques.
We estimate the density $p_t(y)$ of the diffusion at a point $y \in D$
performing an integration by
parts that involves the contribution of the Brownian noise only on an
arbitrarily
small time interval $[t-\delta, t]$. This allows us to gain a free
parameter $\delta$ that we can eventually optimize, and the
appropriate choice of $\delta$ proves to be a key point in our argument.
We do not study here the regularity with respect to initial condition
(which may be the subject of future work).

Our study is motivated by applications to Finance,
in particular by the study of models for financial securities
which rely on SDEs with non-Lipschitz coefficients.
As it is well known, a celebrated process with square-root diffusion
coefficient was proposed
by Cox, Ingersoll and Ross in \cite{CIR} as a model for short interest
rates and was later employed
by Heston in \cite{Hest} to model the stochastic volatilities of assets.
The stochastic-$\alpha\beta\rho$ or SABR model in \cite{sabr} is
based on the
following mixing of local and stochastic volatility dynamics:
\[
\cases{
d X_t
= \sigma_t X_t^{\beta} \,dW^1_t,
\cr
d \sigma_t
= \nu\sigma_t \,dB_t, &\quad $\sigma_0 = \alpha$,}
\]
where $0 < \beta< 1$, $ B_t = \rho W^1_t + \sqrt{1 - \rho^2} W^2_t$,
$(W^1, W^2)$ is a standard Brownian
motion, and $\rho\in[-1, 1]$ is the correlation parameter.
In this paper, we apply our results to one-dimensional SDEs of the form
%
%
\begin{equation}\label{e:localCIR1}
X_t = x + \int_0^t \bigl( a (X_s) - b (X_s) X_s \bigr) \,ds + \int_0^t \gamma
(X_s) X_s^{\alpha} \,dW_t,
\end{equation}
where $\alpha\in[1/2, 1)$ and $a$, $b$ and $\gamma$ are $C^{\infty
}_b$ functions.
When the coefficients~$a$, $b, \sigma$ are constant, the solutions to
this class of equations include
the classical CIR process ($\alpha= 1/2$) and a subclass of the CEV
local volatility diffusions
(when $a = 0$ and $b = - r$).
As pointed out by Bossy and Diop in \cite{BD}, SDEs with square-root terms
and coefficients depending on the level of the state variable arise
as well in the modeling of
turbulent flows in fluid mechanics.
It is well known that for a~CIR process, the density of $X_t$ is known
explicitly.
The main contribution of our results lies in the fact that they apply
to the more general framework of SDEs whose coefficients are functions
of the state variable, thus when explicit computations are no
longer possible.
Theo\-rem~\ref{t:localDensSubLin1} directly applies to (\ref{e:localCIR1})
and allows one to show that~%
$X_t$ admits a smooth density on $(0, +\infty)$: under some additional
conditions on the coefficients (mainly dealing with their asymptotic behavior
at $\infty$ and zero), we give exponential-type upper bounds for the
density at infinity (Proposition~\ref{p:asInfty}) and
study the explosive behavior of the density at zero (Proposition~\ref
{p:asZero}).
The explicit expression of the density for the classical CIR process
shows that our estimates are in the
good range.

The paper is organized follows: in Section \ref{s:generalSmooth} we
present our
main results on SDEs with locally smooth coefficients (Section \ref
{s:mainResults}),
and we collect all the technical elements we need to give their proofs.
In particular, in Section \ref{s:MC} we recall the basic tools of
Malliavin calculus on the Wiener space,
which will be used in Section \ref{s:techPart} to obtain some explicit
estimates of the $L^2$-norms
of the weights
involved in the integration by parts formula.
This is done following some standard techniques of estimation of
Sobolev norms and inverse moments
of the determinant of the Malliavin matrix (as in \cite{Nual} and
\cite{Eul1}, Section 4), but in
our computations we explicitly pop out the dependence with respect to
the coefficients of
the SDE and their derivatives.
This further allows us to obtain the explicit asymptotic estimates on the density.
Section \ref{s:proofs} is devoted to the proof of the theorems stated
in Section~\ref{s:mainResults}.
We employ the Fourier transform argument and the optimized integration
by parts
we have discussed above.
Finally, in Section \ref{s:localCIR} we apply our results to the
solutions of (\ref{e:localCIR1}).

\section{Smoothness and tail estimates of densities for SDEs with
locally smooth coefficients}
\label{s:generalSmooth}

\subsection{Main results} \label{s:mainResults}

In what follows, $b$ and $\sigma_j$ are measurable functions from $\R
^m$ into $\R^m$,
$j = 1,\ldots, d$.
For $y_0 \in\R^m$ and $R > 0$, we denote by $B_R(y_0)$ [resp.,
$\overline{B}_R(y_0)$]
the open (resp., closed) ball $B_R(y_0) = \{ y \in\R^m\dvtx| y - y_0 | <
R \}$
[resp., $\overline{B}_R(y_0) = \{ y \in\R^m\dvtx| y - y_0 | \le R \}$],
where \mbox{$| \cdot|$} stands for the Euclidean norm.
We follow the usual notation denoting $C^{\infty}_b(A)$ the class of
infinitely differentiable
functions on the open set $A \subseteq\R^m$ which are bounded
together with their partial derivatives
of any order.
For a multi-index $\alpha\in\{1,\ldots,m\}^k$, $k \geq1$, $\partial
_{\alpha}$ denotes
the partial derivative $\frac{\partial^{k}}{\partial_{x_{\alpha_1}}
\cdots\partial_
{x_{\alpha_k}}}$.

Let $0 < R \le1 $ and $y_0 \in\R^m$ be given.
We consider the SDE
%
%
\begin{eqnarray}\label{e:base1}
X^i_t
= x^i + \int_0^t b^i(X_s) \,ds
+ \sum_{j=1}^d \int_0^t \sigma^i_j(X_s) \,dW^j_s,\nonumber\\[-8pt]\\[-8pt]
&&\eqntext{t \in[0,T],
i = 1,\ldots, m,}
\end{eqnarray}
for a finite $T > 0$ and $x \in\R^m$, and assume that the following hold:

\begin{enumerate}[(H1)]
\item[(H1)] (\textit{local smoothness})
$b,\sigma_j \in C^{\infty}_b( B_{5 R} (y_0) ; \R^m)$;
\item[(H2)] (\textit{local ellipticity})
$\sigma\sigma^*(y) \geq c_{y_0, R}
I_m$ for every $y \in B_{3 R}(y_0)$, for some $0 < c_{y_0, R} < 1$;
\item[(H3)] existence of strong solutions holds for the
couple $(b, \sigma)$.
\end{enumerate}

Let then $(X_t; t \in[0,T])$ denote a strong solution of (\ref{e:base1}).
Our first main result follows.
\begin{theor} \label{t:localDens1}
Assume \textup{(H1)}, \textup{(H2)} and \textup{(H3)}. Then for any
initial condition $x \in\R^m$ and any $0 < t \le T$, the random vector
$X_t$ admits an infinitely differentiable density $p_{t, y_0}$ on
$B_{R}(y_0)$. Furthermore, for any integer $k \ge3$ there exists a
positive constant $\Lambda_{k}$ depending also on $y_0, R, T$, $m, d$
and on the coefficients of (\ref{e:base1}) such that, setting
\[
P_t (y) =
\Prob\bigl( \inf\{ | X_s - y |\dvtx s \in[(t - 1) \vee t/2, t] \} \leq
3 R \bigr),
\]
then one has
%
%
\begin{equation}\label{e:densEstimate1}
p_{t, y_0} (y) \leq P_t (y_0) \biggl( 1 + \frac{1}{ t^{ m 3 /2 } }
\biggr) \Lambda_{3}
\end{equation}
for any $y \in\Bun$.
Analogously, for every $\alpha\in\{ 1,\ldots, m \}^k$, $k \geq1$,\vspace*{-1pt}
%
%
\begin{equation} \label{e:densDerivEstimate1}
| \partial_{\alpha} p_{t, y_0} (y ) | \leq P_t (y_0) \biggl( 1 + \frac{1}{
t^{ m (2 k + 3) / 2 } } \biggr) \Lambda_{ 2 k + 3}\vspace*{-1pt}
\end{equation}
for every $y \in\Bun$.\vspace*{-2pt}
\end{theor}

The functional dependence of $\Lambda_k$ with respect to $y_0$, $R$,
$T$ and to the bounds on the
coefficients $b$ and $\sigma$ is known explicitly.
We provide the expression of $\Lambda_k$ in Section \ref{s:proofs} in
a more detailed version of Theorem
\ref{t:localDens1} (Theorem \ref{t:localDens}) which we do not give
here for the simplicity of notation.

When the coefficients of (\ref{e:base1}) are smooth outside a compact
ball and have polynomial growth together with their derivatives
therein, according to Theorem \ref{t:localDens1} a smooth density exists outside
the same compact set, and one can deduce some more easily-read bounds on
the tails. More precisely, we consider the following assumptions:\vspace*{-1pt}

\begin{enumerate}[(H1$'$)]
\item[(H1$'$)] There exist $\eta\geq0$ such that $b, \sigma
_j$ are of class $C^{\infty}$
on $\R^m \setminus\overline{B}_{\eta}(0)$, 
and (H2) holds for any $R > 0$ and $y_0$ such that
$\overline{B}_{3 R}(y_0) \subset\R^{m} \setminus\overline{B}_{\eta}(0)$;
\item[(H4)] there exist $q, \overline{q} >
0$, and positive constants
$0 < C_0 <1$ and $C_k$, $k \geq1$, such that for any $\alpha\in\{
1,\ldots, m \}^k$\vspace*{-1pt}
%
%
\begin{equation}\label{e:polGrowth}
| \partial_{\alpha} b^i(y) | + | \partial_{\alpha} \sigma^i_j (y) |
\leq
C_k ( 1 + | y |^q )\vspace*{-2pt}
\end{equation}
and\vspace*{-2pt}
%
%
\begin{equation}\label{e:polElliptic}
\sigma\sigma^*(y)
\geq
C_0 | y |^{ - \overline{q} } I_m\vspace*{-2pt}
\end{equation}
hold for $| y | > \eta$.\vspace*{-2pt}
\end{enumerate}
\begin{theor} \label{t:localDensSubLin1}
Assume \textup{(H1$'$)} and \textup{(H3)}.
\begin{enumerate}[(a)]
\item[(a)] For any initial condition $x \in\R^m$ and for any
$0 < t \le T$,
$X_t$ admits a~smooth density on $\R^{m} \setminus\overline{B}_{\eta}(0)$.
\item[(b)] Assume \textup{(H4)} as well.
Then estimates (\ref{e:densEstimate1}) and (\ref{e:densDerivEstimate1})
hold with \mbox{$R \!=\! 1$ and}\vspace*{-1pt}
%
%
\begin{equation}\label{e:LambdaBound1}
\Lambda_k = \Lambda_k ( y_0 )
:= C_{k, T}
\bigl( 1 + | y_0 |^{ q'_k ( q ) } \bigr),\vspace*{-1pt}
\end{equation}
for every $| y_0 | > \eta+ 5$.
The value of the exponent $q'_k ( q )$ is explicitly known (and
provided in Theorem
\ref{t:localDensSubLin}).
\item[(c)]
If moreover $\sup_{ 0 \le s \le t} | X_s |$ has finite moments of all orders,
then for every \mbox{$p > 0$} and $k \geq1$ there exist positive constants
$C_{k,p, T}$
such that\vspace*{-2pt}
%
%
\begin{eqnarray}\label{e:localDensSubLin1}
| p_t( y ) |
&\leq&
C_{3,p, T}
\biggl( 1 + \frac{ 1 } { t^{ m 3 / 2 } } \biggr)
| y |^{-p},
\nonumber\\[-11pt]\\[-11pt]
| \partial_{\alpha} p_t( y ) |
&\leq&
C_{k,p,T}
\biggl( 1 + \frac{ 1 } { t^{ m ( 2 k + 3 )/2 } } \biggr)
| y |^{-p},\qquad \alpha\in\{1,\ldots, m \}^k,
\nonumber\vspace*{-2pt}
\end{eqnarray}
for every $0 < t \le T$ and every $| y | > \eta+ 5$.\vspace*{-2pt}
\end{enumerate}
In the above, the $C_{k,p,T}$ are positive constants depending on
$k,p,T$ and also on $m, d$ and on the bounds (\ref{e:polGrowth}) and
(\ref{e:polElliptic}) on the coefficients.
\end{theor}

The proofs of these results will be given in Section \ref{s:proofs}.

\subsubsection*{Notation}
Through the rest of the paper, $\langle \cdot, \cdot\rangle $
will denote the Euclidean scalar product in $\R^m$, while the notation
$| \cdot|$ will be\vadjust{\goodbreak} used both
for the absolute value of real numbers and for the Euclidean norm in
$\R^m$.
Furthermore, when $\Theta= \theta_1,\ldots, \theta_{\nu}$ is a
family of parameters,
unless differently specified
by $C_{\Theta}$, we denote a constant depending on the $\theta_i$'s
but not on any of the
other existing variables.
All constants of such a type may vary from line to line, but always
depend only on\vspace*{1pt} the $\theta_i$'s.
For functions of one variable, the $k$th derivative will be denoted by
$\frac{d^k f}{d x^k}$
or $f^{(k)}$.\vspace*{2pt}
We will follow the convention of summation over repeated indexes,
wherever present.

\subsection{Elements of Malliavin calculus} \label{s:MC}

We recall hereafter some elements of Malliavin calculus on the Wiener
space, following \cite{Nual}.

Let $W = (W^1_t,\ldots, W^d_t; t \geq0)$ be a $d$-dimensional
Brownian motion defined on the
canonical space $(\Omega, \mathcal{F}, \Prob)$.
For fixed $T > 0$, let $\Hh$ be the Hilbert space $\Hh= L^2([0,T]; \R^d)$.
For any $h \in\Hh$ we set $W(h) = \sum_{j=1}^d \int_0^T h^j(s) \,dW^j_s$,
and consider the family $\mathcal{S} \subset L^2(\Omega, \mathcal
{F}, \Prob)$ of smooth random
variables defined by
\[
\mathcal{S} = \{ F \dvtx F = f(W(h_1),\ldots, W(h_n));
h_1,\ldots,h_n \in\Hh; f \in
C^{\infty}_{\mathrm{pol}}(\R^n); n \geq1 \},
\]
where $C^{\infty}_{\mathrm{pol}}$ denotes the class of $C^{\infty}$
functions which have polynomial
growth together with their derivatives of any order.

The Malliavin derivative of $F \in\mathcal{S}$ is the $d$-dimensional
stochastic process
$D F = (D^1_r F,\ldots, D^d_r F; r \in[0,T])$ defined by
\[
D^j_r F = \sum_{i=1}^n \frac{\partial f}{\partial x_i}(W(h_1),\ldots,
W(h_n)) h_i^j(r),\qquad j = 1,\ldots, d.
\]
For any positive integer $k$, the $k$th order derivative of $F$ is
obtained by iterating the derivative
operator: for any multi-index $\alpha= (\alpha_1,\ldots, \alpha_k)
\in\{1,\ldots, d\}^k$ and
$(r_1,\ldots, r_k) \in[0,T]^k$, we set
$D^{\alpha_1,\ldots, \alpha_k}_{r_1, \ldots, r_k} F := D_{r_1}^{\alpha
_1} \cdots D_{r_k}^{\alpha_k} F$.
Given $p \geq1$ and positive integer $k$, for every $F \in\mathcal
{S}$ we define the seminorm
\[
\|F\|_{k,p} =
\Biggl( \mathbb{E}[|F|^p]
+ \sum_{h=1}^k \mathbb{E}\bigl[\bigl\|D^{(h)}F\bigr\|_{\Hh^{\otimes
^h}}^p\bigr] \Biggr)^{1/p},
\]
where
\[
\bigl\|D^{(k)}F\bigr\|_{\Hh^{\otimes^k}} =
\biggl( \sum_{|\alpha|=k} \int_{[0,T]^k} |
D^{\alpha_1,\ldots, \alpha_k}_{r_1, \ldots, r_k} F |^2\,
dr_1 \cdots dr_k \biggr)^{1/2},
\]
and the sum is taken over all the multi-indexes $\alpha= (\alpha
_1,\ldots,\alpha_k)
\in\{1,\ldots, d\}^k$.
We denote with $\mathbb{D}^{k,p}$ the completion of $\mathcal{S}$
with respect to the seminorms
$\|\cdot\|_{k,p}$, and we set $ \mathbb{D}^{\infty} = \bigcap_{p \geq
1} \bigcap_{k \geq1} \mathbb{D}^{k,p}$.
We may occasionally refer to $\| F \|_{k,p}$ as the stochastic Sobolev
norm of $F$.

In a similar way, for any separable Hilbert space $V$ we can define the
analogous spaces $\mathbb{D}^
{k,p}(V)$ and $\mathbb{D}^{\infty}(V)$ of $V$-valued random variables
with the corresponding
$\| \cdot\|_{k,p,V}$ semi-norms (the smooth functionals being now of
the form $F = \sum_{j=1}^n F_j v_j$,
where $F_j \in\mathcal{S}$ and $v_j \in V$).
In particular, for any $\R^d$-valued process $(u_s; s \leq t)$ such that
$u_s \in\mathbb{D}^{k,p}$ for all $s \in[0,t]$ and
\[
\|u\|_{\Hh}
+ \sum_{h = 1}^k \bigl\|D^{(h)}u\bigr\|_{\Hh^{\otimes^{h+1}}}
< \infty,\qquad \Prob\mbox{-a.s.},
\]
we have
\[
\|u\|_{k,p,\Hh} =
\Biggl( \mathbb{E}[ \|u\|_{\Hh}^p ]
+ \sum_{h = 1}^k \mathbb{E} \bigl[ \bigl\|D^{(h)} u \bigr\|_{\Hh^{\otimes
^{h+1}}}^p \bigr] \Biggr)^{1/p}.
\]
Finally, we denote by $\delta$ the adjoint operator of $D$.

One of the main applications of Malliavin calculus consists of showing
that the law of a nondegenerate random vector
$F = (F^1,\ldots, F^m) \in(\mathbb{D}^{\infty})^{m}$ admits an
infinitely differentiable density.
The property of nondegeneracy, understood in the sense of the Malliavin
covariance matrix, is introduced in the following:
\begin{defin}
A random vector $F = (F^1,\ldots, F^m) \in(\mathbb{D}^{\infty})^m$,
$m \geq1$, is said to be
nondegenerate if its Malliavin covariance matrix $\sigma_F$, defined by
\[
(\sigma_F)_{i,j} = \langle D F^i, D F^j\rangle_{\Hh},\qquad i,j = 1,\ldots, m,
\]
is invertible a.s. and moreover
\[
\esp[\det(\sigma_F)^{-p}] < \infty
\]
for all $p \geq1$.
\end{defin}

The key tool to prove smoothness of the density for a nondegenerate
random vector is the following
integration by parts formula (cf. \cite{Nual}).
\begin{prop} \label{p:integrByParts}
Let $F = (F^1,\ldots, F^m) \in(\mathbb{D}^{\infty})^m$, $m \geq1$,
be a nondegenerate
random vector.
Let $G \in\mathbb{D}^{\infty}$ and $\phi\in C^{\infty}_{\mathrm
{pol}}(\R^m)$.
Then for any $k \geq1$ and any multi-index $\alpha= (\alpha_1,\ldots,
\alpha_k)
\in\{1,\ldots, m\}^k$ there exists a random variable $H_{\alpha
}(F,G) \in\mathbb{D}^{\infty}$ such that
%
%
\begin{equation}\label{e:integrByParts}
\esp[ \partial_{\alpha} \phi(F) G ] = \esp[\phi
(F) H_{\alpha}(F,G) ],
\end{equation}
where the $H_{\alpha}(F,G)$ are recursively defined by
\begin{eqnarray*}
H_{\alpha}(F,G) &=& H_{(\alpha_k)}\bigl( F, H_{(\alpha_1,\ldots, \alpha
_{k-1})} (F,G)\bigr),
\\
H_{(i)}(F,G) &=& \sum_{j = 1}^m \delta( G (\sigma_F^{-1})_{i,j}
D F^j ).
\end{eqnarray*}
\end{prop}

\subsection{Explicit bounds on integration by parts formula for
diffusion processes} \label{s:techPart}

The notation of this section is somehow cumbersome, as we try to keep
our bounds as general and as
accurate as possible.
The framework will nevertheless considerably simplify in Section \ref
{s:proofs}, when we will give
the proofs of the results stated in Section \ref{s:mainResults}.

Throughout this section, $X = (X_t; t \geq0)$ will denote the unique
strong solution of the SDE\vspace*{-1pt}
%
%
\begin{equation}\label{e:diffusion}
X^i_t =
x^i + \int_0^t B^i(X_s)\,ds + \sum_{j=1}^d \int_0^t
A^i_j(X_s)\,dW^j_s,\qquad
t \geq0, i = 1,\ldots, m,\hspace*{-30pt}\vspace*{-1pt}
\end{equation}
where $x \in\R^m$ and $B^i, A^i_j \in\mathcal{C}^{\infty}_b (\R
^m)$ for
all $i = 1,\ldots, m$ and $j = 1,\ldots, d$.
We assume that the diffusion matrix $A$ satisfies the following
ellipticity condition at
starting point $x$:

(E) $A(x)A(x)^* \geq c_* I_m$,
for some $c_* > 0$, where $\cdot^*$ stays for matrix transposition.
Without loss of generality, we will suppose $c_* < 1$.

We recall that the first-variation process of $X$ is the matrix-valued process\vspace*{-1pt}
\[
(Y_t)_{i,j}
= \frac{\partial X^i_t}{\partial x_j}, \qquad i, j = 1,\ldots, m,\vspace*{-1pt}
\]
which satisfies the following equation, written in matrix form:\vspace*{-1pt}
\[
dY_t =
I_m + \int_0^t \partial B(X_s) Y_s \,ds + \sum_{l=1}^d \int_0^t
\partial A_l(X_s) Y_s \,dW^l_s,\vspace*{-1pt}
\]
where $\partial B$ and $\partial A_l$ are, respectively, the $m \times m$
matrices of components
$(\partial B)_{i,j} = \partial_j B^i$ and $(\partial A_l)_{i,j} =
\partial_j A^i_l$.
By means of It\^o's formula, one shows that $Y_t$ is invertible and
that the inverse $Z_t := Y_t^{-1}$
satisfies the equation\vspace*{-1pt}
%
%
\begin{eqnarray}\label{E:Z}
Z_t& =&
I_m - \int_0^t Z_s \Biggl\{ \partial B(X_s) - \sum_{l=1}^d (\partial
A_l(X_s))^2 \Biggr\} \,ds\nonumber\\[-10pt]\\[-10pt]
&&{} - \sum_{l=1}^d \int_0^t Z_s\, \partial A_l(X_s)\,dW^l_s.\nonumber\vspace*{-1pt}
\end{eqnarray}

\subsubsection*{Additional notation}
For $k \geq0$, we define\vspace*{-1pt}
%
%
\begin{eqnarray}\label{e:coeffBounds}
\nB_k &=&
1 + \sum_{i=1}^m \sum_{0\leq|\alpha|\leq k} \sup_{x \in\R^m}
|\partial_{\alpha}
B^i(x)|,
\nonumber\\[-10pt]\\[-10pt]
\nA_k &=&
1 +
\sum_{i,j}
\sum_{0 \leq|\alpha| \leq k} \sup_{x \in\R^m} |\partial_{\alpha
} A^i_j(x)|,\nonumber\vspace*{-1pt}
\end{eqnarray}
where $| \alpha|$ is the length of the multi-index $\alpha$.
Then, for $p \ge1$ and $t \geq0$ we~set
%
%
\begin{equation}\label{e:ep}
e_p(t) := e^{ t^{p/2} ( t^{1/2} \nB_1 + \nA_1 )^p }
\end{equation}
and\vspace*{2pt}
%
%
\begin{equation}\label{e:epZ}
e^Z_p(t) := e^{ t^{p/2} ( t^{1/2}(\nB_1 + \nA_1^2) + \nA_1
)^p }.\vspace*{2pt}
\end{equation}
The constants in (\ref{e:ep}) and (\ref{e:epZ}) naturally arise when
estimating the moments of the
random variables $X_t$, $Y_t$ and $Z_t$.
Indeed, the results given in the following proposition can be easily
obtained from (\ref{e:diffusion})
and (\ref{E:Z}) applying Burkholder's inequality and Gronwall's lemma.\vspace*{1pt}
\begin{prop}
For every $p > 1$ there exists a positive constant $C_{p, m}$ depending
on $p$ and $m$ but not on
the bounds on $B$ and $A$ and their derivatives such that, for every $0
\leq s \leq t \le T$,\vspace*{2pt}
%
%
\begin{eqnarray} \label{E:XBound}\qquad
\mbox{\hspace*{16.4pt}\textup{(i)\quad}}
\esp
\Bigl[ \sup_{s \leq r \leq t} | X^i_r - X^i_s|^p \Bigr] &\le&
C_{p, m} ( t - s)^{p/2}
\bigl( (t-s)^{1/2} \nB_0 + \nA_0 \bigr)^p,
\\[3pt]
\label{E:ZBound}
\mbox{\textup{(ii)\hspace*{16.4pt}\quad}}
\sup_{s \leq t}
\esp[ |(Z_s)_{i,j}|^p ]
&\leq&
C_{p, m}
e^Z_p(t)^{ C_{p, m} }\vspace*{2pt}
\end{eqnarray}
for all $i,j = 1,\ldots, m$.
\end{prop}

For any $t >0$, the iterated Malliavin derivative of $X_t$ is the
solution of a~linear SDE.
The coefficients of this equation are bounded, and hence it is once
again a straightforward
application of Gronwall's lemma to show that the random variables
$D^{\alpha_1,\ldots, \alpha_k}
_{r_1,\ldots, r_k} X_t$ have moments of any order which are finite and
uniformly bounded in
$r_1,\ldots, r_k$.
This is indeed the content of \cite{Nual}, Theorem 2.2.2.
The following lemma highlights the explicit constants appearing in the
estimates of the $L^p$-norms
of the iterated derivative, expressing them in terms of the bounds
(\ref{e:coeffBounds}) on $A$ and $B$.\vspace*{1pt}
\begin{lem} \label{l:boundDeriv}
For every $k \geq1$ and every $p>1$ there exist a positive integer
$\gamma_{k,p}$ and a positive constant
$C_{k,p}$ depending on $k,p$ but not on the bounds on $B$ and $A$ and
their derivatives such that, for
any $t > 0$,\vspace*{2pt}
%
%
\begin{eqnarray}\label{e:boundDeriv}
&&
\sup_{r_1,\ldots,r_k \leq t} \mathbb{E} [ | \DX
|^p ]\nonumber\\[-6pt]\\[-6pt]
&&\qquad\leq
C_{k,p} \nA_{k-1}^{kp}
( t^{1/2} \nB_k + \nA_k )^{ (k+1)^2 p}
e_p(t)^{ \gamma_{k,p} },\nonumber\vspace*{2pt}
\end{eqnarray}
for all $i=1,\ldots,m$ and $(j_1,\ldots, j_k) \in\{1,\ldots,d\}^k$.
\end{lem}

The proof of this result is based on some standard but rather
cumbersome computations; hence
we leave it for Appendix \ref{a:app1}.
We rather give hereafter the proof of some estimates which follow
easily from Lemma \ref{l:boundDeriv}
and will be useful in the following sections.
\begin{cor} \label{c:boundIntegrDeriv}
For any $k \geq1$ and $p>1$, there exists a positive constant~%
$C_{k,p}$ depending only on $k$ and $p$ such that, for any
$t >0$,
%
%
\begin{eqnarray} \label{e:boundIntegrDeriv}\hspace*{40pt}
\mbox{\textup{(i)}\quad}
\mathbb{E} \bigl[ \bigl\| D^{(k)} X^i_t \bigr\|_{\Hh^{\otimes^k}}^p\bigr]^{1/p}
&\leq& C_{k,p}
t^{k/2}
\nA_{k-1}^{k}
( t^{1/2} \nB_k + \nA_k )^{ (k + 1)^2}\nonumber\\[-8pt]\\[-8pt]
&&\hspace*{0pt}{}\times
e_p(t)^{ \gamma_{k, p} };\nonumber
\\[2pt]
\label{e:boundChainRule}
\mbox{\textup{(ii)}\quad\hspace*{38.5pt}}
\| \phi( X_t ) \|_{k,p}
&\leq& C_{k,p}
| \phi|_k
\bigl( 1 + ( t \vee t^k )^{1/2} \bigr)\nonumber\\[-8pt]\\[-8pt]
&&\hspace*{0pt}{}\times
\nA_{k-1}^{k}
( t^{1/2} \nB_k + \nA_k )^{ ( k + 2 )^2 }
e_p(t)^{ k \gamma_{k, p} },\nonumber
\end{eqnarray}
where \textup{(i)} holds for $i = 1,\ldots, m$ and \textup{(ii)} for any $\phi\in
C^{\infty}(\R^m)$.
\end{cor}
\begin{pf}
(i) Employing the definition of \mbox{$\|\cdot\|_{\Hh^{\otimes^k}}$}
and Lemma \ref{l:boundDeriv}, a simple
computation holds.
\begin{eqnarray*}
&&
\mathbb{E} \bigl[ \bigl\| D^{(k)} X^i_t \bigr\|_{\Hh^{\otimes^k}}^p \bigr]^{1/p}
\\[2pt]
&&\qquad\leq C_{k,p}
\biggl\{
t^{ k ({p}/{2}-1 ) }
\int_{[0,t]^k}
\esp
\Bigl[
\sup_{|\alpha|=k} | D^{\alpha_1,\ldots, \alpha_k}_{r_1,
\ldots, r_k} X^i_t |^p
\Bigr]
\,dr_1 \cdots dr_k
\biggr\}^{1/p}
\\[2pt]
&&\qquad\leq C_{k,p}
t^{k/2}
\nA_{k-1}^{k}
( t^{1/2} \nB_k + \nA_k )^{ (k + 1)^2 }
e_p(t)^{ \gamma_{k, p} / p },
\end{eqnarray*}
hence we get bound (\ref{e:boundIntegrDeriv}).

(ii) We start from the definition of $\| \cdot\|_{k,p}$ and write
%
%
\begin{eqnarray} \label{e:phi1}
\| \phi( X_t ) \|_{k,p}
&=&
\Biggl( \mathbb{E}[ |\phi( X_t )|^p ]
+ \sum_{h=1}^k \mathbb{E} \bigl[ \bigl\| D^{(h)} \phi( X_t ) \bigr\|_{\Hh
^{\otimes^h}}^p \bigr]
\Biggr)^{1/p}
\nonumber\\[-8pt]\\[-8pt]
&\leq&
\| \phi\|_0
+ \sum_{h=1}^k \mathbb{E} \bigl[ \bigl\| D^{(h)} \phi( X_t ) \bigr\|_{\Hh
^{\otimes^h}}^p \bigr]^{1/p}.
\nonumber
\end{eqnarray}
Using the notation introduced in the proof Lemma \ref{l:boundDeriv},
we have
\[
D^{(h)} \phi( X_t )
=
\sum_{ I_1,\ldots, I_{\nu} = \{ 1,\ldots, h \}}
\partial_{k_1} \cdots\partial_{k_{\nu}}
\phi( X_t )
\prod_{l = 1}^{\nu}
D^{ ( \mathrm{card}(I_l) ) } X_t^{k_l},
\]
where, with a slight abuse of notation, we have now written $D^{(h)}$
for the generic derivative
of order $h$.
Repeatedly applying H\"older's inequality for Sobolev norms and using
bound (\ref{e:boundIntegrDeriv}),
we get
\begin{eqnarray*}
&&
\esp\bigl[ \bigl\| D^{(h)} \phi( X_t ) \bigr\|_{\Hh^{{\otimes}^h}}^p
\bigr]\\
&&\qquad\leq
c_{h,p} \mathop{\sum_{h_1,\ldots, h_{\nu} = 1,\ldots, h
}}_{ h_1 + \cdots+ h_{\nu} = h } \esp\Biggl[ \Biggl\|
\partial_{k_1} \cdots\partial_{k_{\nu}} \phi( X_t )
\prod_{l = 1}^{\nu}
D^{ (h_l) } X_t^{k_l}
\Biggr\|_{\Hh^{{\otimes}^h}}^p
\Biggr]
\\[2pt]
&&\qquad\leq
c_{h,p} \| \phi\|_h^p \mathop{\sum_{h_1,\ldots, h_{\nu} =
1,\ldots, h }}_{ h_1 + \cdots+ h_{\nu} = h } \sup_{ i = 1,\ldots, m }
\prod_{l = 1}^{\nu} \esp\bigl[ \bigl\| D^{ (h_l) } X_t^i
\bigr\|_{\Hh^{{\otimes}^{h_l}}}^{2^lp} \bigr]^{1/2^l}
\\[-1pt]
&&\qquad\leq
c_{h,p} \| \phi\|_h^p \bigl\{ t^{h/2} \nA_{h-1}^{h} ( t^{1/2}
\nB_h + \nA_h )^{ ( h + 2)^2 } \bigr\}^{p} e_p(t)^{h \gamma_{h, p}}.
\end{eqnarray*}
By means of this bound, from (\ref{e:phi1}) we get the desired
estimate when setting
$C_{k,p} \geq\max\{ c_{h,p} \dvtx h \leq k \}$.\vspace*{-2pt}
\end{pf}

We need a last preliminary result on the inverse moments of the
determinant of the Malliavin
covariance matrix of $X_t$.
This result is again achieved with some standard arguments, but,
as in Lemma \ref{l:boundDeriv}, the next lemma finds out the explicit
constants appearing
in the estimate of the $L^p$-norms of $\det(\sigma_{X_t})^{-1}$.\vspace*{-2pt}
\begin{lem} \label{l:covMatrix}
For every $p>1$ and $t > 0$,
%
%
\begin{equation}\label{e:covMatrix}
\mathbb{E} [ | {\det\sigma_{X_t}} |^{-p} ]^{1/p} \leq C_{p, m, d}
e^Z_{4(mp+1)}(t)^{ C_{p, m, d} } K_m(t, c_*),
\end{equation}
where\vspace*{-3pt}
\[
K_m(t, c_*) = 1 + \biggl( \frac{4}{t c_*} + 1 \biggr)^m + \frac{1}{
c_*^{2(m + 1) } } ( t^{1/2}\|B\|_0\|A\|_2^3 + \|A\|_1^2 )^{ 2(m+1) }
\]
for some positive constant $C_{p, m, d}$ depending on $p, m$ and $d$
but not on the bounds on
$B$ and $A$ and their derivatives.\vspace*{-2pt}
\end{lem}

The proof is once again postponed to Appendix \ref{a:app1}.

We now come to the main result of this section.
We give an estimate of the $L^2$-norm of the random variables
$H_{\alpha}$
involved in the integration by parts formula (\ref{e:integrByParts}),
when $F = X_t$.
The proof follows the arguments of \cite{Eul1}, proof of Lemma 4.11,
but is given in the general setting of an integration by parts
of order $k \in\mathbb{N}$, and moreover it takes advantage of the
explicit bounds
which have been obtained in Corollary \ref{c:boundIntegrDeriv} and
Lemma \ref{l:covMatrix}.

$\!\!$We give this result employing some slightly more compact notation,~\mbox{defining}
\begin{eqnarray*}
P_k (t) &=& t^{1/2} | B |_k + | A |_k,
\\
P^{A}_k(t) &=& | A |_k P_{k+1} (t).\vspace*{-2pt}
\end{eqnarray*}

\begin{theor} \label{p:Hbound}
For every $k \geq1$ there exists a positive constant $C_k = C_{k, m, d}$
such that, for any multi-index $\alpha\in\{1,\ldots,m\}^k$, any $G
\in\mathbb{D}^{\infty}$
and $t > 0$,
%
%
\begin{eqnarray}\label{e:Hbound}
&&
\| H_{\alpha}(X_t,G) \|_{0,2}\nonumber\\[-1pt]
&&\qquad\leq C_k
\|G\|_{k,2^{k+1}}
\bigl( t^{ -{k}/{2} } \vee t^{ {k (k - 1)}/{2} } \bigr)
( t^m K_m(t, c_*) )^{{k (k + 3)}/{2}}
\\[-1pt]
&&\qquad\quad{} \times
( P^{A}_k )^{ \phi_k }
\bigl( e_8(t) \vee e_{ 2^{k + 2} }(t) \bigr)^{ C_k }
\bigl( e^Z_{ 32m + 4 }(t) \vee e^Z_{2^{ k + 4 } m + 4}(t) \bigr)^{
C_k },
\nonumber
\end{eqnarray}
where $K_m(t, c_*)$ has been defined in Lemma \ref{l:covMatrix}, and\vspace*{-1pt}
\[
\phi_k = 3 m (k + 4)^2.\vspace*{-1pt}
\]
\end{theor}

\begin{rem}
Estimate (\ref{e:Hbound}) is rather involved.
For our purposes, the most important elements are the dependence with
respect to time of the factor $t^{ -{k}/{2} } \vee t^{ {k (k
- 1)}/{2} }$ and the coefficient $P^{A}_k$ containing the bounds on the
derivatives of the coefficients.
We remark that the factor $t^m K_m(t, c_*)$ is bounded for $t$ close to zero.
Moreover, when $t < 1$, the factor $t^{ -{k}/{2} } \vee t^{
{k (k - 1)}/{2} }$ reduces to $t^{-{k}/{2}}$.\vspace*{-5pt}
\end{rem}

\begin{pf*}{Proof of Theorem \ref{p:Hbound}}
We write $\sigma_t = \sigma_{X_t}$ for simplicity of notation.
We first use the continuity of $\delta$ (see \cite{Eul1}, Proposition
4.5) and H\"older's
inequalities for Sobolev norms 
to obtain
%
%
\begin{eqnarray} \label{e:H1}\quad
&&\bigl\| H_{(\alpha_1,\ldots,\alpha_{k})}(X_t,G) \bigr\|_{0,2}\nonumber\\[-2pt]
&&\qquad= \bigl\| H_{(\alpha_k)}\bigl(X_t, H_{(\alpha_1,\ldots,\alpha_{k-1})}(X_t,G)\bigr) \bigr\|_{0,2}
\nonumber\\[-10pt]\\[-10pt]
&&\qquad= \Biggl\| \sum_{j=1}^m
\delta\bigl( H_{(\alpha_1,\ldots,\alpha_{k-1})}(X_t,G)
(\sigma_t^{-1})_{\alpha_k,j} D X^j_t \bigr)
\Biggr\|_{0,2}
\nonumber\\[-4pt]
&&\qquad\leq C_m \bigl\| H_{(\alpha_1,\ldots,\alpha_{k-1})}(X_t,G) \bigr\|_{1,4}
\sum_{j=1}^m \| (\sigma_t^{-1})_{\alpha_k,j} \|_{1,8} \| D X^j_t
\|_{1,8,\Hh}.\nonumber\vspace*{-2pt}
\end{eqnarray}
To estimate the last factor we can directly use the definition of $\|
\cdot\|_{k,p,\Hh}$
and apply Corollary \ref{c:boundIntegrDeriv}.
The major part of the efforts in the rest of the proof will be targeted
on the estimation
of $\| (\sigma_t^{-1})_{i,j} \|_{k, p}$.

We claim that for any $k \geq1$, $p > 1$ and for all $i, j = 1,\ldots, m$,
%
%
\begin{eqnarray}\label{e:invMMatrix}
\| (\sigma_t^{-1})_{i,j} \|_{k,p}
&\leq& c_{k, p}
( t^{-1} \vee t^{{k}/{2}-1} )
( t^m K_m(t, c_*) )^{1 + k}
\nonumber\\[-10pt]\\[-10pt]
&&{}
\times P^A_k ( t )^{ \phi'_k + 2 (k + 4)^2 }
e_p(t)^{ c_{k, p} } e_{4(mp+1)}^Z(t)^{ c_{k, p} },
\nonumber\vspace*{-1pt}
\end{eqnarray}
where\vspace*{-5pt}
\[
\phi'_k = 2 (k + 1) (m - 1),
\]
and $c_{k, p}$ is a positive constant depending also on $m, d$ but not
on $t$ and on the bounds
on $B$ and $A$ and their derivatives.
Iterating process (\ref{e:H1}) and repeatedly using estimates (\ref
{e:invMMatrix})
and (\ref{e:boundIntegrDeriv}), one easily obtains the desired estimate\vspace*{-1pt}
\begin{eqnarray*}
&&
\bigl\| H_{ (\alpha_1,\ldots,\alpha_{k}) }( X_t, G ) \bigr\|_{0,2}\\
&&\qquad\leq C_{k, m, d}
\| G \|_{k,2^{k+1}}
( t^m K_m(t, c_*) )^k
\\[-3pt]
&&\qquad\quad{} \times\prod_{h=1}^k
( t^{-1} \vee t^{{h/2}-1} ) (t \vee t^h
)^{{1/2}}
( t^m K_m(t, c_*) )^h
\\[-3pt]
&&\hspace*{59pt}{} \times
P^A_k ( t )^{ k ( \phi'_k + 2 (k + 4)^2 + (k + 1)^2 ) }
\\[-3pt]
&&\qquad\quad{} \times
\prod_{h=1}^k
e_{ 2^{h + 2} }(t)^{ c_{h,m,d} }
e^Z_{ 4 (2^{h + 2} m + 1 ) } (t)^{ c_{h,m,d} }
\\[-3pt]
&&\qquad\leq C_{k, m, d}
\|G\|_{k,2^{k+1}}
( t^m K_m(t, c_*) )^{{k (k + 3)}/{2}}
\\[-1pt]
&&\qquad\quad{} \times
\bigl( t^{-{k}/{2}} \vee t^{{k (k - 1)}/{2}} \bigr)
P^{A}_k(t)^{ \phi_k }
\\[-1pt]
&&\qquad\quad{} \times
\bigl( e_8(t) \vee e_{ 2^{k + 2} }(t) \bigr)^{ C_{k, m, d} }
\bigl( e^Z_{ 32m + 4 }(t) \vee e^Z_{ 2^{ k + 4 } m + 4}(t) \bigr)^{
C_{k, m, d} }.
\end{eqnarray*}

\textsc{Proof of (\ref{e:invMMatrix}).}\quad
We follow \cite{Eul1}, proof of Lemma 4.11. We start from the
definition of $\| \cdot\|_{k,p}$
and write\vspace*{-2pt}
%
%
\begin{equation}\label{e:defNormInvMMatrix}\qquad
\| (\sigma_t^{-1})_{i,j} \|_{k,p}
= \Biggl( \mathbb{E}[ | (\sigma_t^{-1})_{i,j} |^p ]
+ \sum_{h=1}^k \mathbb{E}\bigl[ \bigl\| D^{(h)}(\sigma_t^{-1})_{i,j}
\bigr\|_{\Hh^{\otimes^h}}^p \bigr] \Biggr)^{1/p}.\vspace*{-2pt}
\end{equation}
For the first term, we simply use Cramer's formula for matrix inversion,
\[
| (\sigma_t^{-1})_{i,j} | = ( \det\sigma_t)^{-1} \sigma_t^{(i,j)},
\]
where $\sigma_t^{(i,j)}$ denotes the $(i, j)$ minor of $\sigma_t$.
We then apply H\"older's inequality and bounds
(\ref{e:boundIntegrDeriv}) and (\ref{e:covMatrix}) and get
%
%
\begin{eqnarray}\label{e:boundInvMMatrix0}
\esp[ | (\sigma_t^{-1})_{i,j} |^p ]
&\leq& c^{(1)}_{p,m}
\bigl\{
\esp[ \det( \sigma_t )^{-2p} ]
\esp\bigl[ \bigl| \sigma_t^{(i,j)} \bigr|^{-2p} \bigr]
\bigr\}^{1/2}
\nonumber\\
&\leq& c^{(1)}_{p,m}
\Bigl\{
\esp[ \det( \sigma_t )^{-2p} ] \esp\Bigl[ \sup_i
\| D X^i_t \|_{\Hh}^{4(m-1)p} \Bigr]
\Bigr\}^{1/2}
\nonumber\\
&\leq& c^{(1)}_{p, m, d}
t^{-p} ( t^m K_m(t, c_*) )^p
\\
&&{} \times\{ \nA_0 ( t^{1/2} \nB_1 + \nA_1
)^4 \}^{ 2 (m - 1) p}
\nonumber\\
&&{} \times e_p(t)^{ c^{(1)}_{p, m, d} } e^Z_{4(mp+1)}(t)^{
c^{(1)}_{p, m, d} },\nonumber
\end{eqnarray}
where $K_m(t, c_*)$ is the constant defined in Lemma \ref{l:covMatrix}.
To estimate the second term, as done in
\cite{Eul1}, proof of Lemma 4.11, we iterate the chain rule for $D$
\[
D (\sigma_t^{-1})_{i,j} = - \sum_{a,b=1}^m (\sigma_t^{-1})_{i,a}
D(\sigma_t)_{a,b} (\sigma_t^{-1})_{b,j}.
\]
We take advantage of the notation introduced in the proof of Lemma \ref
{l:boundDeriv}
and for
$(\beta_1,\ldots,\beta_k) \in\{1,\ldots,m\}^k$, $k \geq1$,
we write
%
%
\begin{eqnarray} \label{e:derMMatrix}\qquad\quad
&&
| D^{\beta_1,\ldots,\beta_k}_{r_1,\ldots,r_k} (\sigma_t^{-1})_{i,j}|
\nonumber\\
&&\qquad\leq
\sum_{I_1 \cup\cdots\cup I_{\nu} = \{1,\ldots,k\} }
\mathop{\sum_{a_1,\ldots,a_{\nu}=1}}_{b_1,\ldots,b_{\nu}=1}^m
| (\sigma_t^{-1})_{i,a_1} (\sigma_t^{-1})_{b1,a_2}
\cdots\nonumber\\[-8pt]\\[-8pt]
&&\qquad\quad\hspace*{109pt}{}\times(\sigma
_t^{-1})_{b_{\nu-1},a_{\nu}} (\sigma_t^{-1})_{b_{\nu},j} |
\nonumber\\
&&\qquad\quad\hspace*{106.7pt}{}\times
\bigl| D^{\beta(I_1)}_{r(I_1)} (\sigma_t)_{a_1,b_1} \cdots D^{\beta
(I_{\nu})}_{r(I_{\nu})} (\sigma_t)_{a_{\nu},b_{\nu}}
\bigr|.\nonumber
\end{eqnarray}
We repeatedly apply H\"older's inequality for Sobolev norms to (\ref
{e:derMMatrix}) and get
%
%
\begin{eqnarray}\label{e:derInvMMatrix}
&&
\esp\bigl[ \bigl\| D^{(k)} (\sigma_t^{-1})_{i,j} \bigr\|_{\Hh^{\otimes^k}}^p\bigr]\nonumber\\
&&\qquad\leq
c^{(2)}_{k,p,m}
\mathop{\sum_{k_1,\ldots, k_{\nu} = 1,\ldots, k}}_{k_1 + \cdots+
k_{\nu} = k}
\Bigl\{
\mathop{\sup_{a, a_1,\ldots,a_{\nu} = 1,\ldots, m}}_{b, b_1,\ldots
,b_{\nu} = 1,\ldots, m}
\esp\bigl[ \bigl\|
(\sigma_t^{-1})_{a,b}^{\nu+ 1}
\nonumber\\
&&\qquad\quad\hspace*{170.7pt}{}\times
D^{(k_1)} (\sigma_t)_{a_1,b_1} \cdots\nonumber\\[-2pt]
&&\qquad\quad\hspace*{182.3pt}{}\times D^{(k_{\nu})}
(\sigma_t)_{a_{\nu},b_{\nu}}
\bigr\|_{\Hh^{\otimes^k}}^p\bigr]
\Bigr\}
\\[-1pt]
&&\qquad\leq c^{(2)}_{k,p,m}\nonumber\\[-4pt]
&&\qquad\quad{}\times
\mathop{\sup_{k_1,\ldots, k_{\nu}=1,\ldots,k}}_{k_1 + \cdots+
k_{\nu} = k}
\Biggl\{
\sup_{a,b=1,\ldots,m} \esp\bigl[ | (\sigma_t^{-1})_{a,b} |^{(\nu
+1) p} \bigr]
\nonumber\\[-4pt]
&&\hspace*{148pt}{}\times
\prod_{l=1}^{\nu}
\sup_{ a_l, b_l = 1,\ldots, m }
\esp\bigl[ \bigl\| D^{(k_l)} (\sigma_t)_{a_l,b_l} \bigr\|_{\Hh
^{\otimes^{k_l}}}^{2^l p} \bigr]^{1/2^l}
\Biggr\},
\nonumber
\end{eqnarray}
where, as in the proof of Corollary \ref{c:boundIntegrDeriv}, we have
written $D^{(k_l)}$ for
the generic derivative of order $k_l$.
To estimate $D^{(k_l)} (\sigma_t)_{a_l, b_l}$ we use bound (\ref
{e:boundIntegrDeriv}) and get
%
%
\begin{eqnarray}\label{e:boundMMatrix}\qquad
&&\esp\bigl[ \bigl\| D^{(k)} (\sigma_t)_{i,j} \bigr\|_{\Hh^{\otimes^k}}^p
\bigr]\nonumber\\[-1pt]
&&\qquad\leq
\esp\Biggl[ \Biggl\| \sum_{h=0}^k
\pmatrix{k \cr h}
\int_0^t D^{(h)} D_s X^i_t \cdot D^{(k-h)} D_s X^j_t
\Biggr\|_{\Hh^{\otimes^k}}^p
\Biggr]
\nonumber\\[-1pt]
&&\qquad\leq c^{(3)}_{k,p}
\sum_{h=0}^k
\esp\bigl[ \bigl\| D^{(h)} D X^i_t\bigr\|^{2p}_{\Hh^{\otimes^{h+1}}}\bigr]^{1/2}
\\[-1pt]
&&\qquad\quad\hspace*{33pt}{} \times
\esp\bigl[ \bigl\| D^{(k-h)} D X^j_t\bigr\|^{2p}_{\Hh^{\otimes
^{k-h+1}}}\bigr]^{1/2}
\nonumber\\[-1pt]
&&\qquad\leq c^{(3)}_{k,p}
t^{( {k}/{2} + 1 ) p}
\nA_k^{ (k + 2) p}
( t^{1/2} \nB_{k+1} + \nA_{k+1} )^{ 2 (k + 2)^2 p}
e_p(t)^{2 \gamma_{k+1, p}},\nonumber
\end{eqnarray}
where we have once again applied H\"older's inequality for Sobolev norms
in the second
step.

Using (\ref{e:boundMMatrix}) together with (\ref{e:derInvMMatrix}),
bound (\ref{e:boundInvMMatrix0})
and (\ref{e:defNormInvMMatrix}) and observing that $t^m K_m(t)$ is
greater than one for all
the values of $t$, we finally obtain
\begin{eqnarray*}
\| (\sigma_t^{-1})_{i,j} \|_{k,p}
&\leq&
c_{k, p}
( t^{-1} \vee t^{{k/2}-1} )
( t^m K_m(t, c_*) )^{1 + k}
\\[-1pt]
&&{} \times\nA_k^{ \phi'_k + k (k + 2) }
( t^{1/2} \nB_{k+1} + \nA_{k+1} )^{ \phi'_k + 2 (k + 4)^2 }
\\[-1pt]
&&{} \times e_p(t)^{ c_{k, p} } e_{4(mp+1)}^Z(t)^{ c_{k, p} },
\end{eqnarray*}
for a positive constant $c_{k, p}$ depending also on $m, d$.
Estimate (\ref{e:invMMatrix}) follows.
\end{pf*}

\subsection{\texorpdfstring{Proofs of Theorems \protect\ref{t:localDens1} and
\protect\ref{t:localDensSubLin1}}{Proofs of Theorems 2.1 and 2.2}}
\label{s:proofs}

We now come to the proof of the results stated in Section
\ref{s:mainResults}. We recall that an $\R^m$-valued random vector $X$
is said to admit a density on an open set $A \in\R^m$ if $\mathcal{L}_X
|_A$ possesses a~density, $\mathcal{L}_X$ being the law of~$X$. It is
equivalent to say that
%
%
\begin{equation}\label{e:locDensDef}
\esp[ f(X) ] = \int_{\R} f(x) p(x) \,dx
\end{equation}
holds for all $f \in C_b(\R)$ such that $\operatorname{supp}(f) \subset A$, for some
positive $p \in L^1(A)$.

We refer to the setting of Section \ref{s:mainResults}.
We recall that $X = (X_t; t \in[0, T])$ denotes a strong solution of
%
%
\begin{eqnarray}\label{e:base}
X^i_t
= x^i + \int_0^t b^i(X_s) \,ds
+ \sum_{j=1}^d \int_0^t \sigma^i_j(X_s) \,dW^j_s,\nonumber\\[-8pt]\\[-8pt]
&&\eqntext{t \in[0,T], i = 1,\ldots, m,}
\end{eqnarray}
where $b$ and $\sigma$ satisfy the assumptions \textup{(H1)}--\textup{(H3)}.
For $k \geq1$ and $f \in C^{k}(\R^m)$, we denote
%
%
\begin{equation}\label{e:coeffNorms}
| f |_{k, B_R(y_0)} =
1
+ \sum_{|\alpha|\leq k} \sup_{x \in B_R(y_0)} |\partial_{\alpha}
f(x)|,
\end{equation}
where the sum is taken over all the multi-index $\alpha\in\{1,\ldots, m
\}^k$.
Let us define the following ``local'' version of the constants
appearing in the
estimates of the previous section:
\begin{eqnarray*}
P_{k}(t,y_0) &=& t^{1/2} | b |_{k, B_{5 R}(y_0)} + | \sigma|_{k, B_{5
R}(y_0)},\\
P^{\sigma}_{k}(t,y_0) &=& | \sigma|_{k, B_{5 R}(y_0)} P_{k+1}(t,y_0),
\\
P^{Z}_{1}(t,y_0) &=& t^{1/2} \bigl( | b |_{1, B_{5 R}(y_0)} + | \sigma
|^2_{1, B_{5 R}(y_0)} \bigr)
+ | \sigma|_{1, B_{5 R}(y_0)},
\\
P_m^C(t,y_0)
&=& \bigl( t^{1/2} | b |_{0,B_{5 R}(y_0)} | \sigma|_{2,B_{5 R}(y_0)}^3
+ | \sigma|_{1,B_{5 R}(y_0)}^2 \bigr)^{2(m+1)},
\\
C_m(t,y_0) &=&
t^m
+ \frac{4^m}{ c_{y_0}^{2 (m + 1)} }
\bigl(1 + P^C_m(t,y_0)\bigr),
\\
e_p(t, y_0)
&=& \exp( t^{p/2} P_1(t, y_0)^p ),\\
e_p^Z(t, y_0)
&=& \exp( t^{p/2} P_1^Z(t, y_0)^p ).
\end{eqnarray*}
In order to prove Theorem \ref{t:localDens1}, we simplify this rather
heavy notation introducing a
constant that contains the factors appearing in estimate (\ref
{e:Hbound}) in Theorem \ref{p:Hbound}
(recall the constant $\phi_k$ defined there)
\begin{eqnarray*}
\Theta_k(t, y_0, \gamma)
&=&
C_m(t, y_0)^{ {mk (mk + 3)}/{2} }
P^{\sigma}_{mk} (t,y_0)^{ \phi_{ m k } + (m k + 2)^2 }
\\
&&{}
\times
\bigl( e_8(t, y_0) \vee e_{ 2^{mk + 2} }(t, y_0) \bigr)^{ \gamma}\\
&&{}
\times
\bigl( e^Z_{ 32m + 4 }(t, y_0) \vee e^Z_{ 2^{ m k + 4 } m + 4}(t, y_0)
\bigr)^{ \gamma}.
\end{eqnarray*}

As addressed in Section \ref{s:mainResults}, the following theorem is
a more detailed version of Theorem \ref{t:localDens1}.
In particular, it provides the explicit expression of the constant
$\Lambda_k$ appearing in estimates (\ref{e:densEstimate1}) and (\ref
{e:densDerivEstimate1}).
\begin{theor}\label{t:localDens}
Assume \textup{(H1)}, \textup{(H2)} and \textup{(H3)}.
Then, for any initial condition $x \in R^{m}$ and any $0 < t \le T$,
the random vector $X_t$ admits an infinitely differentiable density
$p_{t, y_0}$ on $B_{R}(y_0)$.
Furthermore, for every $k \geq1$
there exists a positive constant $C_k = C_{k, m, d}$ such that,
setting
%
%
\begin{equation}\label{e:LambdaDef}
\Lambda_k( t, y_0 )
=
C_k
R^{ -m k }
\bigl(
P_0( t, y_0 )^{ m k }
+ \Theta_{ k }( t, y_0 , C_k )
\bigr)
\end{equation}
and
\[
P_t (y) =
\Prob\bigl( \inf\{ | X_s - y |\dvtx s \in[(t - 1) \vee t/2, t] \} \leq
3 R \bigr),
\]
then one has
%
%
\begin{equation}\label{e:densEstimate}
p_{t, y_0} (y)
\leq P_t (y_0)
\biggl( 1 + \frac{1}{ t^{ m 3 /2 } } \biggr)
\Lambda_3 ( t \wedge1, y_0 )
\end{equation}
for every $y \in\Bun$.
Analogously, for any $\alpha\in\{ 1,\ldots, m \}^k$, $k \ge1$,
%
%
\begin{equation} \label{e:densDerivEstimate}
| \partial_{\alpha} p_{t, y_0} (y ) |
\leq P_t (y_0)
\biggl( 1 + \frac{1}{ t^{ m 3 /2 } } \biggr)
\Lambda_{ 2 k + 3 }( t \wedge1, y_0 )
\end{equation}
for every $y \in\Bun$.
\end{theor}

To prove this result we rely on the following classical criterion for
smoothness of laws based on a
Fourier transform argument (cf. \cite{Nual}, Lemma 2.1.5).
\begin{prop} \label{p:FourInv}
Let $\mu$ be a probability law on $\R^m$, and $\widehat{\mu}(\xi)
= \int_{\R} e^{i \langle \xi, y\rangle}\times\break \mu(dy)$
its characteristic function.
If $\widehat{\mu}$ is integrable, then $\mu$ is absolutely
continuous w.r.t. the Lebesgue measure,
and
%
%
\begin{equation}\label{e:fInv}
p(y)
= \frac{1}{(2\pi)^m} \int_{\R^m} e^{-i\langle \xi, y\rangle} \widehat{\mu
}(\xi) \,d \xi
\end{equation}
is a continuous version of its density.
If moreover
%
%
\begin{equation}\label{e:fourierDecay}
\int_{\R^m} | \xi|^k |\widehat{\mu}(\xi)| \,d \xi< \infty
\end{equation}
holds for any $k \in\mathbb{N}$, then $p$ is of class $C^{\infty}$
and for any multi-index $\alpha=
(\alpha_1,\ldots, \alpha_k) \in\{1,\ldots, m \}^k$,
\[
\partial_{\alpha} p (y)
= (-i)^k
\int_{\R}
\Biggl( \prod_{j=1}^k \xi^{\alpha_j} \Biggr)
e^{-i\langle \xi, y\rangle} \widehat{\mu}(\xi) \,d \xi.
\]
\end{prop}
\begin{pf*}{Proof of Theorem \ref{t:localDens}}
\textit{Step} 1 (``\textit{localized}'' \textit{characteristic function}).
Fix a $t$ in $(0,T]$. Let $\phi_R \in C^{\infty}_b(\R^m)$ be such
that $1_{ B_{R}(0) }
\leq\phi_R \leq1_{ B_{2 R}(0) }$ and $ | \phi_R |_k \leq2^{k} R^{-k}$.
We first observe that if $ m_0 = \esp[ \phi_R ( X_t - y_0) ]$ is
zero, then it just follows that
$p \equiv0$ is a density for $X_t$ on $B_{R} (y_0)$.
Otherwise, we consider~$\mathcal{L}_{t,y_0}$ the law on $\R^m$ such that
%
%
\begin{equation}\label{e:condLaw}
\int_{\R^m} f(y) \mathcal{L}_{t, y_0}(dy)
= \frac{1}{m_0}
\esp[ f(X_t) \phi_R ( X_t - y_0 ) ],
\end{equation}
for all $f \in C_b(\R^m)$.
If $\mathcal{L}_{t,y_0}$ possesses a density, say $\pt'$, it follows that
$\pt(y) := m_0 \pt'$
is a density for $X_t$ on $B_R(y_0)$.
Indeed, for any $f \in C_b$ such that $\operatorname{supp}(f) \subset B_R (y_0)$,
(\ref{e:condLaw}) implies
\begin{eqnarray*}
\int_{\R^m} f(y) \pt(y) \,dy
&=&
\int_{\R^m} f(y) m_0 \pt'(y) \,dy
\\
&=&
m_0 \int_{\R^m} f(y) \mathcal{L}_{t, y_0}(dy)
\\
&=&
\esp[ f(X_t) ].
\end{eqnarray*}
If the characteristic function of $\mathcal{L}_{t,y_0}$
\[
\pcap(\xi)
= \int_{\R^m} e^{i \langle\xi, y\rangle} \mathcal{L}_{t,y_0}(dy)
= \frac{1}{m_0}
\esp\bigl[ e^{ i\langle \xi, X_t\rangle} \phi_R ( X_t - y_0 ) \bigr]
\]
is integrable, then by Proposition \ref{p:FourInv} $\mathcal
{L}_{t,y_0}$ admits a density.
Hence, we focus on the integrability of $\pcap$; in particular, we
show that condition (\ref{e:fourierDecay})
of Proposition~\ref{p:FourInv} holds true for all $k \in\mathbb{N}$.

Moreover, the inversion formula (\ref{e:fInv}) yields the
representation for $\pt$
%
%
\begin{eqnarray}\label{e:locDens}
\pt(y)
:\!&=& m_0 \pt'(y)
= \frac{m_0}{(2 \pi)^m} \int_{\R^m} e^{-i\langle \xi,y\rangle }
\pcap(\xi) \,d \xi
\nonumber\\[-8pt]\\[-8pt]
&=& \frac{1}{(2 \pi)^m} \int_{\R^m} e^{-i\langle \xi,y\rangle }
\esp\bigl[ e^{i \langle \xi,X_t\rangle} \phi_R ( X_t - y_0 ) \bigr] \,d \xi.
\nonumber
\end{eqnarray}

\textit{Step} 2 (\textit{localization}). We define the coefficients
%
%
\begin{eqnarray}\label{e:truncCoeff}
\overline{b}{}^i(y)
&=& b^i\bigl(\psi(y - y_0)\bigr),
\nonumber\\[-8pt]\\[-8pt]
\overline{\sigma}{}^i_j(y)
&=& \sigma^i_j\bigl(\psi(y - y_0)\bigr),
\nonumber
\end{eqnarray}
where $\psi\in C^{\infty}(\R^m; \R^m)$ (a \textit{truncation}
function) is defined by
\[
\psi(y)
= \cases{
y, &\quad if $| y | \leq4 R$,\vspace*{2pt}\cr
5 \dfrac{y}{| y |}, &\quad if $| y | \geq5 R$,}
\]
and $\psi(y) \in\overline{B}_{5 R}(0)$ for all $y \in\R^m$.
$\psi$ can be defined in such a way that, for all $i = 1,\ldots, m$,
$\| \psi^i \|_1
\leq1$ and $\| \psi^i \|_k \leq2^{k-2} R^{- (k - 1)}$ for all $k
\geq2$.
As a~consequence of~(H1),
the $\overline{b}$ and $\overline{\sigma}$ defined in this way are
$C^{\infty}_b$-extensions of $b |_{B_{4 R} (y_0)}$
and $\sigma|_{B_{4 R} (y_0)}$.
Furthermore, there exist constants $c = c_{k,m}$ such that
%
%
\begin{eqnarray}\label{e:truncCoeffBounds}
| \overline{b}{}^i |_k
&\leq&
c_{k,m} R^{- (k - 1)}
| b^i |_{k, B_{5 R} (y_0)},
\nonumber\\[-8pt]\\[-8pt]
| \overline{\sigma}{}^i_j |_k
&\leq&
c_{k,m} R^{- (k - 1)}
| \sigma^i_j | _{k, B_{5 R} (y_0)}
\nonumber
\end{eqnarray}
and by (H2), for any $y \in\Bdeux$ the matrix $\overline{\sigma}(y)$
is elliptic
%
%
\begin{equation}\label{e:truncMatrixEllipt}
\overline{\sigma}\,
\overline{\sigma}^*(y)
\geq c_{y_0, R} I_m,\qquad y \in\Bdeux.
\end{equation}
For $y \in\R^m$ we denote by $\overline{X}(y) = (\overline{X}_s(y)
; 0 \leq s \leq t)$ the unique strong
solution of the equation
%
%
\begin{eqnarray}\label{e:localizedX}
\overline{X}{}^i_s(y) =
y^i +
\int_0^s \overline{b}{}^i(\overline{X}_u(y)) \,du
+ \sum_{j=1}^d \int_0^s \overline{\sigma}{}^i_j(\overline{X}_u(y))
\,dW^j_u,\nonumber\\[-8pt]\\[-8pt]
&&\eqntext{0 \leq s \leq t,
i = 1,\ldots,m.}
\end{eqnarray}

Let now $0 < \delta< t/2 \wedge1$.
We employ an up-down crossing argument to estimate the increments of $X$
in the neighborhood of $y_0$
by replacing them with the increments of $\Xbar$.
More precisely, let
$\nu= \nut$ and $\tau= \taut$ be the stopping times defined by
%
%
\begin{eqnarray}\label{e:stopTimes}
\nut
&=& \inf\{ s \geq t - \delta\dvtx X_s \in B_{3 R} (y_0) \},
\nonumber\\[-8pt]\\[-8pt]
\taut
&=& \inf\{ s \geq\nut\dvtx X_s \notin B_{4 R} (y_0) \}
\nonumber
\end{eqnarray}
and $\inf\{ \varnothing\} = \infty$.
Suppose that $\phi_R ( X_t - y_0 ) > 0$, so that $X_t \in B_{ 2 R} (
y_0 )$
and $\nu< t$.
On this set, if $\nu> t - \delta$, then $| X_{ t \wedge\tau} -
X_{\nu} | \ge R$.
This implies\break $| \Xbar_{ t \wedge\tau- \nu} ( X_{\nu} )
- X_{\nu}| = | X_{ t \wedge\tau} - X_{\nu} | \ge R$.
Here we are employing the fact that on the interval $[\nu, \tau]$,
$X$~stays in $B_{4 R}(y_0)$, hence in the region where the truncated
coefficients
$\overline{b}, \overline{\sigma}$ coincide with the original ones
$b, \sigma$.
On this interval, both $X$ and $\overline{X}$ satisfy (\ref
{e:localizedX}) for which pathwise
uniqueness holds; hence we can replace $X$ by $\Xbar$ and
employ the flow property for $\overline{X}$.
Notice that flow property may not hold true for $X$ [due to possible
lack of uniqueness for the
couple $(b,\sigma)$], but it always does for $\overline{X}$.

Analogously, if $\nu\!=\! t\! -\! \delta$ and $ \tau\!<\! t $, then $| X_{ \tau
} \!-\! X_{\nu} | \!=\!
| \Xbar_{ \tau- \nu} ( X_{\nu} ) \!-\! X_{\nu}|\! \ge\! R$.
In both cases,
${\sup_{ 0 \le s \le\delta} }| \Xbar_s ( X_{\nu} ) - X_{\nu}| \ge R$.
Hence, we conclude that
\begin{eqnarray*}
\{ \phi_R( X_t - y_0 ) > 0 \}
&=&
\{ \phi_R( X_t - y_0 ) > 0, t - \delta= \nu< t < \tau\}
\\
&&{}
\cup
\Bigl\{
\phi_R( X_t - y_0 ) > 0,
{\sup_{ 0 \le s \le\delta}}
| \Xbar_s ( X_{\nu} ) - X_{\nu} | \ge R
\Bigr\}
\end{eqnarray*}
and $\pcap$ rewrites as
\begin{eqnarray*}
m_0 \pcap(\xi)
&=&
\esp
\bigl[
e^{i\langle \xi, X_t\rangle } \phi_R ( X_t - y_0 )
1_{ \{ \phi_R ( X_t - y_0 ) > 0, \sup_{0 \le s \le\delta} |
\Xbar_s ( X_{\nu} ) - X_{\nu} | \ge R \} }
\bigr]
\\
&&{}
+
\esp
\bigl[
e^{i\langle \xi, X_t\rangle } \phi_R ( X_t - y_0 )
1_{ \{ \phi_R ( X_t - y_0 ) > 0, t - \delta= \nu< t < \tau\} }
\bigr].
\end{eqnarray*}
We now claim that for all $q > 0$ the following estimate holds:
%
%
\begin{eqnarray} \label{e:stopTimeEstimates}
&&\Prob
\Bigl(
\phi_R( X_t - y_0 ) > 0,
{\sup_{ 0 \le s \le\delta}}
| \Xbar_s ( X_{\nu} ) - X_{\nu}| \ge R
\Bigr)
\nonumber\\[-8pt]\\[-8pt]
&&\qquad\le
c_{q, m}
R^{ - q}
\delta^{q/2}
P_0(\delta, y_0)^{q}
\Prob\Bigl( {\inf_{t-\delta\le s \le t} }| X_s - y_0 | \leq3 R
\Bigr),\nonumber
\end{eqnarray}
for some positive constant $c_{q,m}$.
Estimate (\ref{e:stopTimeEstimates}) will be proved later on.

On the other hand,
%
%
\begin{eqnarray} \label{e:CF3}\qquad
&&\esp
\bigl[
e^{ i\langle \xi, X_t\rangle  } \phi_R ( X_t - y_0 )
1_{ \{ \phi_R ( X_t - y_0 ) > 0, t - \delta= \nu< t < \tau\} }
\bigr]
\nonumber\\
&&\qquad=
\esp
\bigl[
\esp
\bigl[
e^{ i\langle \xi, \Xbar_{\delta}(y) \rangle } \phi_R \bigl( \Xbar_{\delta}(y) - y_0 \bigr)
|
X_{t - \delta} = y
\bigr]
1_{ \{ t - \delta= \nu< t < \tau\} }
\bigr]
\\
&&\qquad\leq
\Prob( | X_{t - \delta} - y_0 | < 3 R )
\sup_{y \in\Bdeux}
\bigl|\esp\bigl[ e^{i\langle \xi, \Xbar_{\delta}(y)\rangle } \phi_R \bigl( \Xbar_{\delta
}(y) - y_0 \bigr) \bigr]
\bigr|.\nonumber
\end{eqnarray}

\textit{Step} 3 (\textit{integration by parts}).
We apply integration by parts formula (\ref{e:integrByParts}) to
estimate the last
term in (\ref{e:CF3}).
By (\ref{e:truncCoeffBounds}), (\ref{e:truncMatrixEllipt}) and Lemma
\ref{l:covMatrix},
$\Xbar_{\delta}(y)$ is a smooth and nondegenerate random vector for
any $\delta>0$ and $y \in\Bdeux$.
Then, for a given $k \geq1$ we define the multi-index
\[
\alpha= ( \underbrace{1,\ldots, 1}_{k\ \mathrm{times}},\ldots,
\underbrace{m ,\ldots, m}_
{k \ \mathrm{times}} ),
\]
such that $| \alpha| = km$.
Hence, recalling that $\partial_{x_k} e^{ i \langle \xi, x \rangle  }
= i \xi^k e^{ i \langle \xi, x \rangle }$,
%
%
\begin{eqnarray}\label{e:intParts}
&&\bigl|
\esp\bigl[ e^{i\langle \xi, \Xbar_{\delta}(y)\rangle } \phi_R \bigl( \Xbar_{\delta
}(y) - y_0 \bigr) \bigr]
\bigr|
\nonumber\\
&&\qquad\leq
\frac{1}{ {\prod_{i=1}^m }| \xi^{i}|^k }
\bigl|\esp\bigl[\partial_{\alpha} e^{i\langle \xi, \Xbar_{\delta}(y)\rangle } \phi_R \bigl( \Xbar
_{\delta}(y) - y_0 \bigr)
\bigr]\bigr|
\\
&&\qquad\leq
\frac{1}{ {\prod_{i=1}^m | \xi^{i}|^k }}
\esp\bigl[ \bigl| H_{\alpha} \bigl( \Xbar_{\delta}(y), \phi_R \bigl( \Xbar
_{\delta}(y) - y_0 \bigr) \bigr) \bigr| \bigr],\nonumber
\end{eqnarray}
for any $y \in B_{3 R} (y_0)$.

We need to separately estimate $\| \phi( \Xbar_{\delta}(y) - y_0
)\|_{| \alpha|, 2^{| \alpha| + 1}}$.
By Corollary \ref{c:boundIntegrDeriv}, this is given by
\begin{eqnarray*}
\bigl\| \phi\bigl( \Xbar_{\delta}(y) - y_0 \bigr)\bigr\|_{mk, 2^{mk+1}}
&\leq&
c_{k, m}
R^{-mk}
( 1 + \delta^{1/2} )
| \sigma|_{ mk -1, B_{5 R} (y_0) }^{ mk }
\\
&&{} \times
P_{ m k } (y_0, \delta)^{ ( m k + 2 )^2}
e_{ 2^{ mk + 1 } } (\delta)^{ c_{k, m} }
\\
&\leq&
c_{k, m}
R^{-mk}
P^{\sigma}_{mk} (y_0, \delta)^{ ( m k + 2 )^2 }
e_{ 2^{ mk + 1 } } (\delta)^{ c_{k, m} }
\end{eqnarray*}
for some positive constant $c_{k, m}$.
Then, from (\ref{e:stopTimeEstimates}), (\ref{e:CF3}), (\ref
{e:intParts}) and Theorem \ref{p:Hbound}
it follows that
%
%
\begin{equation}\label{e:cfEstim}
m_0 | \pcap(\xi) |
\leq
C_{k, q}
P_R(\delta, t, y_0)
I_{k, q} (\xi, \delta, y_0)
\end{equation}
for some constant $C_{k, q}$ depending also on $m$ and $d$, with
\[
P_R(\delta, t, y_0)
=
\Prob\Bigl( {\inf_{ t - \delta\leq s \leq t}} | X_s - y_0 | \leq3 R
\Bigr)
\]
and
\[
I_{k, q} (\xi, \delta, y_0)
=
R^{ - q }
\delta^{ q / 2 }
P_0 (\delta, y_0)^q
+
\frac{ R^{ - m k } }{ \prod_{i=1}^m | \xi^i |^k }
\delta^{- mk/2}
\Theta_k(\delta, y_0, C_{k, q}).
\]
Estimate (\ref{e:cfEstim}) holds simultaneously for any $\xi\in\R
^m$, $0 < \delta< t/2 \wedge1$, \mbox{$q > 0$} and $k \geq1$.
The constant $\Theta_k(\delta, y_0, C_{k, q})$ appears when applying
estima\-te~(\ref{e:Hbound}).

\textit{Step} 4 (\textit{optimization}).
We show that for any $\xi$ and any $l \geq1$, $\delta$ can always be chosen
in such a way that there exist $q$ and $k$ such that $I_{k, q} (\xi,
\delta, y_0)$
goes to zero at $\infty$ faster than $( {\prod_{i=1}^m} | \xi^i
| )^{ - (l + 2) }$.

Denoting $\| \xi\| = {\prod_{i=1}^m} | \xi^i |$, we set
\[
\delta:= \delta(\xi) = t/2 \wedge1 \wedge\| \xi\|^{ -a }
\]
for some $a > 0$ that is to be identified hereafter.
For this choice of $\delta$,
\[
P_R( \delta(\xi) , t, y_0)
\leq
\Prob\Bigl( {\inf_{ t/2 \vee(t - 1) \leq s \leq t }} | X_s - y_0 |
\leq3 R \Bigr) = P_t (y_0)
\]
and
%
%
\begin{eqnarray} \label{e:cfEstim2}\quad
I_{k, q} (\xi, \delta(\xi), y_0)
&\leq&
R^{ - q }
\bigl(
\| \xi\|^{- {q a}/{2}}
\wedge
( t \wedge1)^{ {q}/{2} }
\bigr)
P_0( t \wedge1, y_0 )^q
\nonumber\\
&&{} +
R^{ - m k }
\bigl(
\| \xi\|^{-k (1 - {m a}/{2})}
\vee
\| \xi\|^{-k} (t \wedge1)^{- {m k}/{2}}
\bigr)\\
&&\hspace*{10pt}{}\times
\Theta_k( t \wedge1, y_0, C_{k,q}),\nonumber
\end{eqnarray}
since $\delta\to P_0( \delta, y_0 )$ and $\delta\to\Theta_k(
\delta, y_0, C_{k, q} )$
are increasing; hence $P_0( \delta(\xi), y_0 ) \leq P_0( t \wedge1,
y_0 )$ and the same holds for $\Theta_k$.

We consider the leading terms determining the decay of $I_{k, q} (\xi,
\delta(\xi), y_0)$ with
respect to $\xi$ and impose
%
%
\begin{equation}\label{e:optim}
\frac{q a}{2} = k \biggl(1 - \frac{m a}{2}\biggr).
\end{equation}
Setting $a = 1/m$, (\ref{e:optim}) yields $q = m k$, hence $\frac{q
a}{2} =
k (1 - \frac{m a}{2}) = \frac{k}{2}$.
Therefore, we get the bound
%
%
\begin{eqnarray} \label{e:holdsForAnyK}
I_{k, q^*_k} (\xi, \delta(\xi), y_0)
&\leq&
R^{ - m k } \bigl( \|
\xi\|^{-k/2} \wedge (t \wedge1)^{mk/2} \bigr) P_0(t \wedge1, y_0)^{ m k }
\nonumber\\
&&{}+
R^{ - m k } \bigl( \| \xi\|^{-k/2} \vee (t \wedge1)^{-mk/2} \|
\xi\|^{-k} \bigr) \\
&&\hspace*{10pt}{}\times\Theta_k(t \wedge1, y_0, C_{k, q^*_k})\nonumber
\end{eqnarray}
with $q^*_k = m k$. Estimate (\ref{e:holdsForAnyK}) holds for any $k
\ge1$ and $\xi\neq 0$, and then it proves that the function
$p_{t,y_0}(y)$ defined in (\ref{e:locDens}) is in fact well defined and
infinitely differentiable with respect to $y$.

Let us come to estimate (\ref{e:densEstimate}). We take
(\ref{e:locDens}) and cut off the integration over a region $I$ of
finite Lebesgue measure on which $ \| \xi\| = {\prod_{i=1}^m} | \xi^i
|$ remains smaller than a given constant. That is, we write\vspace*{-1pt}
\begin{eqnarray*}
p_{t, y_0} ( y )
&=&
\frac{ 1 }{ (2 \pi)^m }
\int_{\R^m} e^{-i\langle \xi,y\rangle }
\esp\bigl[ e^{i\langle \xi, X_t\rangle } \phi_R ( X_t - y_0 ) \bigr]
\,d \xi
\\[-1pt]
&\leq&
\frac{ 1 }{ (2 \pi)^m } \biggl[ \int_{ I } \esp[ \phi_R ( X_t - y_0 )] \,d \xi
\\[-1pt]
&&\hspace*{34.6pt}{}
+ \int_{ I^c } e^{-i\langle \xi, y\rangle }
\esp\bigl[ e^{i\langle \xi, X_t\rangle } \phi_R ( X_t - y_0 ) \bigr]
\,d \xi
\biggr]
\\[-1pt]
&\leq&
\frac{ 1 }{ (2 \pi)^m }
\biggl[
\Prob( | X_t - y_0 | < 2 R )
\lambda_m( I )
\\[-1pt]
&&\hspace*{34.6pt}{}
+
C_{k, q^*_k}
P_t (y_0)
\int_{ I^c }
I_{k, q^*_k} (\xi, \delta(\xi), y_0)
\,d \xi
\biggr],
\end{eqnarray*}
where $\lambda_m$ denotes the Lebesgue measure on $\R^m$.
As we have seen, the last term is such that
\begin{eqnarray*}
&&
\int_{ I^c }I_{k, q^*_k} (\xi, \delta(\xi), y_0)
\,d \xi
\\[-1pt]
&&\qquad\le
R^{- m k }
P_0( t \wedge1, y_0 )^{ m k }
\int_{ I^c }
| \xi|^{-k/2}
\,d \xi
\\[-1pt]
&&\qquad\quad{}+ R^{- m k }
\Theta_k(t \wedge1, y_0, C_{k, q^*_k})
\biggl(
(t \wedge1)^{- m k /2}
\int_{ I^c \cap\{ \xi\dvtx| \xi| < (t \wedge1)^{-m} \} }
| \xi|^{-k}
\,d \xi
\\[-1pt]
&&\qquad\quad\hspace*{170pt}{} + \int_{ I^c \cap\{ \xi\dvtx| \xi| \geq(t \wedge1)^{-m} \} }
| \xi|^{-k/2}
\,d \xi
\biggr).
\end{eqnarray*}
Now, since
\[
\int_{ I^c \cap\{ \xi\dvtx| \xi| \geq(t \wedge1)^{-m} \} } | \xi
|^{-k/2} \,d \xi
\leq
\int_{ I^c } | \xi|^{-k/2} \,d \xi
=
c^{(1)}_k < \infty
\]
and
\[
(t \wedge1)^{- m k /2}
\int_{ I^c \cap\{ \xi\dvtx| \xi| < (t \wedge1)^{-m} \} } | \xi|^{-k}
\,d \xi
\leq
(t \wedge1)^{- m k /2} c^{(2)}_k < \infty
\]
hold for any $k \geq3$, we then take $k = 3$ and get the estimate
\[
p_{t, y_0} ( y ) \leq C^* P_t( y_0 ) \bigl[ 1 + R^ {- 3 m } \bigl( P_0(
t \wedge1, y_0 )^{ 3 m } + (t \wedge1)^{- 3 m /2} \Theta_3( t \wedge1,
y_0 , C_{m, d}) \bigr) \bigr]
\]
for every $y \in\Bun$, for a positive constant $C^*$,
estimate (\ref{e:densEstimate}) then follows.
For estimate (\ref{e:densDerivEstimate})
on the derivatives we proceed in the same way, observing that for
$\alpha\in\{1,\ldots, m \}^l$,
$| \xi|^{ - k / 2 } \times\prod_{j=1}^l | \xi^{\alpha_j} |$ is
integrable at $\infty$ as soon
as $k \ge2 l + 3$.

\textsc{Proof of (\ref{e:stopTimeEstimates}).}\quad
We remark that $\{ \phi_R (X_t - y_0) \} \subseteq\{ t-\delta\le\nu
\le t \} \subseteq\{ t-\delta
\le\nu\le t, X_{\nu} \in\overline{B}_{3R} (y_0) \}$, hence\vspace*{-2pt}
\begin{eqnarray*}
&&\Prob \Bigl( \phi_R( X_T - y_0 ) > 0, {\sup_{ 0 \le s \le\delta}} | \Xbar_s (
X_{\nu} ) - X_{\nu}| \ge R \Bigr)
\\[-1pt]
&&\qquad
\le\Prob \Bigl( t-\delta \le\nu\le t, X_{\nu} \in\overline{B}_{3R}
(y_0), {\sup_{ 0 \le s \le\delta}} | \Xbar_s ( X_{\nu} ) - X_{\nu}| \ge R\Bigr)
\\[-1pt]
&&\qquad\le
R^{ -q }
\Prob( t-\delta\le\nu\le t)
\sup_{ y \in\overline{B}_{3R}(y_0) }
\esp\Bigl[
{\sup_{ 0 \leq s \leq\delta}} | \overline{X}_s (y) - y |^q
\Bigr].
\end{eqnarray*}
Using boundedness of coefficients of (\ref{e:localizedX}), it is easy
to show that
\[
\esp \Bigl[ {\sup_{ 0 \leq s \leq\delta} }| \overline{X}_s (y) - y |^q
\Bigr] \leq c_{q, m} \delta^{q/2} P_0 (\delta, y_0)^q
\]
for some positive constant $c_{q, m}$.
\end{pf*}

As addressed in Section \ref{s:mainResults}, the constants appearing
in the definition of $\Lambda_k$ (\ref{e:LambdaDef}) can be
considerably simplified under assumptions (H1$'$) and
(H4), resulting in some polynomial-type bounds.
The following result corresponds to Theorem~\ref{t:localDensSubLin1};
in the presents statement, we explicitly give the expression of the
exponent $q'_k(q)$ appearing in bound (\ref{e:LambdaBound1}).
\begin{theor} \label{t:localDensSubLin}
Assume \textup{(H1$'$)} and \textup{(H3)}.
\begin{enumerate}[(a)]
\item[(a)] For any initial condition $x \in\R^m$ and any $0 < t \le
T$, $X_t$ admits a smooth density on
$\R^{m} \setminus\overline{B}_{\eta}(0)$.
\item[(b)] Assume \textup{(H4)} as well.
Then the constant $\Lambda_k$ defined in Theorem \ref{t:localDens} is
such that
%
%
\begin{equation}\label{e:LambdaBound}
\Lambda_k (t, y_0)
\le C_{k, T}
\bigl( 1 + | y_0 |^{ q'_k ( q ) } \bigr),
\end{equation}
for every $0 < t \le T$ and every $| y | > \eta+ 5$. The exponent
$q'_k ( q )$ is
worth
\[
q'_k ( q ) = m k ( \overline{q} + 4) (m + 1) (m k + 3) + 2 q m \bigl(
\phi_{m k} + (m k + 2)^2 \bigr).
\]
\item[(c)]
If moreover $\sup_{ 0 \le s \le t} | X_s |$ has finite moments of all orders,
then for every \mbox{$p > 0$} and every $k \geq1$ there exist positive
constants $C_{k,p,
T}$ such that
%
%
\begin{eqnarray}\label{e:localDensSubLin}
| p_t( y ) |
&\leq&
C_{3,p, T}
\biggl( 1 + \frac{ 1 } { t^{ m 3 /2 } } \biggr)
| y |^{-p},
\nonumber\\[-10pt]\\[-10pt]
| \partial_{\alpha} p_t( y ) |
&\leq&
C_{k, p, T}
\biggl( 1 + \frac{ 1 } { t^{ m ( 2 k + 3 )/2 } } \biggr)
| y |^{-p},\qquad \alpha\in\{1,\ldots, m \}^k,
\nonumber
\end{eqnarray}
for every $0 < t \le T$ and every $| y | > \eta+ 5$.
\end{enumerate}
The $C_{k,p, T}$ are positive constants depending also on $m, d$ and on
bounds (\ref{e:polGrowth}) and (\ref{e:polElliptic}) on the coefficients.
\end{theor}

\begin{pf}
(a)
We no longer need to distinguish between $y_0$ and the (close) point
$y$ where the density is\vadjust{\goodbreak}
evaluated; hence we just set $y = y_0$ and consider suitable radii.
For $| y | > \eta$, we set $R_{y} = \frac{1}{10} \operatorname{dist} (y, \overline
{B}_{\eta} (0) ) \wedge1$.
By (H1$'$), $b$ and~$\sigma$ are of class $C^{\infty}_b$ on
$B_{5 R_y} (y)$
and satisfy (H2) on $B_{3 R_y} (y)$.
From Theorem~\ref{t:localDens} it follows that $X_t$ admits a smooth
density on $B_{R_y} (y)$.
This holds true for every ball $B_{R_y} (y)$ with center $y$ in $\R^m
\setminus\overline{B}_{\eta}(0)$;
hence statement (a) follows.

(b)
Without loss of generality, we take $R = 1$.
As a consequence of (\ref{e:polGrowth}),
the constants introduced before Theorem \ref{t:localDens}
can be bounded as follows, for $0 \le t \le T$ and $| y | > \eta+
5$:\vspace*{-2pt}
%
%
\begin{eqnarray}\label{e:constantsBounds}
P_{k}(t, y)
&\leq& c^{(1)}_k \bigl( 1 + ( | y | + 5 )^q \bigr)
\leq c^{(1)}_k | y |^q,
\nonumber\\[-1pt]
P^{\sigma}_{k}(t, y) \vee P^{Z}_{1}(t, y)
&\leq&
c^{(1)}_k | y |^{2 q},
\nonumber\\[-1pt]
P_m^C(t, y)
&\leq&
c^{(1)} ( | y |^4 + | y |^2 )^{2 (m + 1)}
\leq
c^{(1)} | y |^{ 8 (m + 1) },
\\[-1pt]
C_m(t, y)
&\leq&
\frac{ c^{(1)} }{ C_0^{2 (m + 1)} }
| y |^{2 \overline{q} (m + 1) }
| y |^{ 8 (m + 1) }\nonumber\\[-1pt]
&\leq&
c^{(1)}
| y |^{2 ( \overline{q} + 4 ) (m + 1) },\nonumber
\end{eqnarray}
for some constants $c^{(1)}$ and $c^{(1)}_k$
depending also on $m, q$ and on the bounds on $b, \sigma$ and their
derivatives in
(\ref{e:polGrowth}) and (\ref{e:polElliptic}).

The exponential factors $e$ and $e^Z$
must be treated on a specific basis.
Indeed, $e_{\cdot}(t, y)$ and $e_{\cdot}^Z(t, y)$ may explode
when $| y | \rightarrow+ \infty$.
Nevertheless, explosion can be avoided stepping further into the
optimization procedure set up in the proof
of Theorem \ref{t:localDens}.
More precisely, we restart from step 4
and force the state variable $y$ to appear in the choice of $\delta$,
setting\vspace*{-1pt}
\[
\delta(\xi, y) = t/2 \wedge1 \wedge|\xi|^{-a} \wedge| y |^{ - 4 q }.
\]
Now, whatever the value of $p$ is, $e_{p}$ and $e_{p}^Z$ are reduced to
\begin{eqnarray*}
e_p( \delta(\xi, y) ) &\le&
e_p( 1 \wedge| y |^{- 4 q}, y )\\[-1pt]
&\leq&
\exp{ \bigl( c^{(1)}_1 ( 1 \wedge| y |^{- 2 q p} )
( 1 + 1 \wedge| y |^{- 2 q} )^p | y |^{ q p } \bigr) }
\\[-1pt]
&\leq&
\exp{ \bigl( c_{p} ( 1 \wedge| y |^{- 2 q p} ) | y |^{ q
p } \bigr) }
\\[-1pt]
&\leq&
\exp{ ( c_{p} ) }\vspace*{-2pt}
\end{eqnarray*}
and\vspace*{-2pt}
\begin{eqnarray*}
e_p^Z( \delta(\xi, y) )
&\le&
e_p^Z(1 \wedge| y |^{-4}, y)\\[-1pt]
&=&
\exp{ \bigl( c^{(1)}_1 ( 1 \wedge| y |^{- 2 q p} )
( 1 + 1 \wedge| y |^{- 2 q} )^p | y |^{2 q p} \bigr)}
\\[-1pt]
&\leq&
\exp{ ( c_p ) }.
\end{eqnarray*}
We then
perform the integration over $\xi$ as done in the last step of the
proof of Theorem \ref{t:localDens}, and
employing (\ref{e:constantsBounds})
we obtain estimate (\ref{e:LambdaBound}) for $\Lambda_k$, for $| y |
> \eta+ 5$.
The value of $q'_k(q)$ is obtained from the definition of $\Theta_k$
and~(\ref{e:constantsBounds}).

(c)
From boundedness of moments of ${\sup_{s \leq t} }| X_s |$, for any
interval $I_t \subseteq[0, t]$
one can easily deduce the estimate
%
%
%
\begin{eqnarray}\label{e:MarkovProba}
&&\Prob( \inf\{ | X_s - y | \dvtx s \in I_t \} \leq3 )
\nonumber\\
&&\qquad\leq
\Prob( \sup\{ | X_s |\dvtx s \in I_t \} \geq| y | - 3 )
\\
&&\qquad\leq
\frac{ 1 }{ ( | y | - 3 )^r }
\esp\Bigl[ \sup_{s \le t} | X_s |^r \Bigr]
\leq
c^{(2)}_r
\frac{ 1 }{ | y |^r },\nonumber
\end{eqnarray}
for any $r > 0$, $ 0 \le t \le T$ and $| y | > 3$.
It is then easy to obtain the desired estimate on $p_t$ with Theorem
\ref{t:localDens}: for a given $p > 0$,
we employ (\ref{e:densEstimate}) with $y_0 = y$ and (\ref
{e:MarkovProba}) with $r > p + q'_{3} ( q )$. 
Similarly, to obtain the estimate on derivatives, one employs (\ref
{e:densDerivEstimate}) and (\ref{e:MarkovProba}) with $r > p + q'_{2 k
+ 3} ( q )$.
\end{pf}

\section{A square root-like (CIR/CEV) process with local coefficients}
\label{s:localCIR}

We apply our results to the solution of 
(\ref{e:localCIR1}).
We will be able to refine the polynomial estimate on the density at $+
\infty$ giving
exponential-type upper bounds.
Under some additional assumptions on the coefficients, we also study
the asymptotic behavior of the density at zero,
that is, the point where the diffusion coefficient is singular.

We first collect some basic facts concerning the solution of (\ref
{e:localCIR1}).
Let us recall the SDE
%
%
\begin{equation}\label{e:localCIR}
\cases{
d X_t = \bigl( a( X_t) - b( X_t ) X_t \bigr)\, dt + \gamma( X_t ) X_t^{\alpha}
\,dW_t, &\quad $t \geq0, \alpha\in[1/2, 1)$,
\cr
X_0 = x \geq0.}\hspace*{-30pt}
\end{equation}
When $\alpha= 1/2$ and $a, b$ and $\gamma$ are constant, the solution
to (\ref{e:localCIR}) is the
celebrated Cox--Ingersoll--Ross process (see \cite{CIR}), appearing in
finance as a~model for short interest rates.
It is well known that, in spite of the lack of globally
Lipschitz-continuous coefficients,
existence and uniqueness of strong solutions hold for the equation of a CIR
process.
If $a \geq0$, the solution stays a.s. in $\R_+ = [0, +\infty)$;
furthermore, a solution starting at $x >0$ stays a.s. in $\R_> = (0,
+\infty)$ if the Feller condition
$2 a \geq\gamma^2$ is achieved (cf. \cite{LL} for details).
The following proposition gives the (straightforward) generalization of
the previous statements to the case
of coefficients $a$, $b$, $\gamma$ that are functions of the
underlying process.
The proof is left to Appendix \ref{a:app2}.
\begin{prop} \label{p:localCIR}
Assume:
\begin{enumerate}[(s0)]
\item[(s0)] $\alpha\in[1/2,1)$; $a, b$ and $\gamma$ $\in C^{1}_b$
with $a(0) \geq0$ and $\gamma(x)^2 > 0$ for every $x > 0$.
\end{enumerate}
Then, for any initial condition $x \geq0$ there exists a unique strong
solution 
to~(\ref{e:localCIR}) which is such that $\Prob( X_t \geq0; t \geq0
) = 1$.
Let then $x > 0$ and $\tau_0 = \inf\{ t \geq0 \dvtx X_t = 0 \}$, with $\inf\{
\varnothing\} = \infty$.
\begin{itemize}
\item If $\alpha> 1/2$ and
\begin{enumerate}[(s1)$'$]
\item[(s1)$'$] $a(0) > 0$ and $z \mapsto\frac{1}{\gamma^2(z)
z^{2 \alpha-1}}$ is integrable
at $0^+$,
\end{enumerate}
then\vspace*{-2pt}
%
%
\begin{equation}\label{e:noZero}
\Prob( \tau_0 = \infty) = 1.\vspace*{-2pt}
\end{equation}
\item If $\alpha= 1/2$ and
\begin{enumerate}[(s2)]
\item[(s1)] $\frac{1}{\gamma^2}$ is integrable at zero,
\item[(s2)] there exists $\overline{x} > 0$ such that
$\frac{ 2a(x) }{ \gamma(x)^2 } \ge1$ for $0 < x < \overline{x}$,\vspace*{-2pt}
\end{enumerate}
then the same conclusion on $\tau_0$ holds.\vspace*{-3pt}
\end{itemize}
\end{prop}

When $X$ is a CIR process, the moment-generating function of $X_t$ can
be computed explicitly, leading to the knowledge of the density.
Setting $L_t = ( 1 - e^{-bt} ) \gamma^2 / 4 b$, then $X_t/L_t$ follows
a noncentral chi-square
law with $\delta= 4 a / \gamma^2$ degrees of freedom and parameter
$\zeta_t = 4 x b /
( \gamma^2 ( e^{b t} - 1 ) )$ (recall that $x$ is here the initial condition).
The density of $X_t$ is then given by (cf.~\cite{LL})\vspace*{-4pt}
\[
p_t( y ) = \frac{ e^{ -\zeta_t/2 } }{ 2^{\delta/2} L_t } e^{ - y / (2
L_t) } \biggl( \frac{y}{L_t} \biggr)^{ \delta/2 - 1 }
\sum_{n=0}^{\infty} \frac{ ( {y}/({4 L_t}))^n }{ n! \Gamma(\delta/2+n)
} \zeta_t^n,\qquad y > 0.\vspace*{-2pt}
\]
We incidentally remark that $p_t$ is in general unbounded, since
$y^{\delta/2-1}$ diverges at zero when $\delta/2-1 = 2a/\gamma^2 - 1$
is negative (in fact, fixed a value of $\delta/2-1$, there
exists a $n \geq0$ such that $\frac{d^n}{d y^n} p_t$ is\vspace*{1pt} unbounded).

The standard techniques of Malliavin calculus cannot be directly
applied to study the existence of
a smooth density for the solution of (\ref{e:localCIR}), as
the diffusion coefficient in general is not (depending on $\gamma$)
globally Lipschitz continuous.
Actually, Alos and Ewald \cite{AE} have shown that if $X$ is CIR
process, then $X_t$, $t > 0$,
belongs to $\mathbb{D}^{1,2}$ when the Feller condition $2 a \geq
\gamma^2$ is achieved.
Higher order of differentiability (in the Malliavin sense) can be
proven, requiring a stronger condition on $a$ and $\gamma$, and the
authors apply these results
to option pricing within the Heston model.
If we are interested in density estimation, the results of the previous
sections allow us to overcome
the problems related to the singular
behavior of the diffusion coefficient and to directly establish the
existence of a smooth density,
independently from any Feller-type condition [provided that (s0) is satisfied].
More precisely, we can give the following preliminary result:\vspace*{-3pt}
\begin{prop}[(Preliminary result)] \label{p:localCIRDens}
Assume \textup{(s0)} and let $a, b, \gamma$ be of class $C^{\infty}_b$.
Let $X = (X_t; t \geq0)$ be the strong solution of (\ref{e:localCIR})
starting at $x \ge0$.
For any $t >0$, $X_t$ admits a smooth density $p_t$ on $(0, +\infty)$.
$p_t$ is such that $\lim_{ y \rightarrow\infty} p_t(y) y^p = 0$ for
any $p > 0$.\vspace*{-2pt}
\end{prop}

\begin{pf}
It is easy so see that, under the current assumptions, the drift and
diffusions coefficients of (\ref{e:localCIR}) satisfy (H1$'$) with
$\eta = 0$ and (H4) with\vadjust{\eject} $q = 1$. (H3)~holds as well, by Proposition
\ref{p:localCIR}. As the coefficients have sub-linear growth, for any
$t > 0$, $\sup_{s \le t} X_s$ has\vspace*{1pt} finite moments of any order. The
conclusion follows from Theorem \ref{t:localDensSubLin}(c).\vspace*{-2pt}
\end{pf}

\subsection{\texorpdfstring{Exponential decay at $\infty$}{Exponential decay at infinity}} \label{s:asInfty}
In order to further develop our study of the density, we could take
advantage of some of the generalized-chaining tools settled by Viens
and Vizcarra in \cite{VV}.
In particular notice that, in order to estimate the density by means of
Theorem \ref{t:localDensSubLin}, we need to deal with the probability
term $P_t(y)$ appearing therein.
For our present purposes, we can rely on alternative strategies
involving time-change arguments and the existence of quadratic
exponential moments for suprema of Brownian motions (Fernique's
theorem) in the current section, and a detailed analysis of negative
moments of the process $X$ in Section
\ref{s:zero}.

From now on, condition (s0) is assumed, the coefficients $a, b, \gamma
$ are of class~$C^{\infty}_b$
and $(X_t; t\le T)$ denotes the unique strong of (\ref{e:localCIR}) on
$[0,T]$, \mbox{$T > 0$}.
We make explicit the dependence with respect to the initial condition
denoting $p_t(x, \cdot)$ the density at
time $t$ of $X$ starting at $x \ge0$.
The following result improves Proposition~\ref{p:localCIRDens} in the
estimate of the density for $y \to\infty$.\vspace*{-2pt}
\begin{prop} \label{p:asInfty}
Assume that\vspace*{-1pt}
%
%
\begin{equation}\label{e:bCondition1}
\varliminf_{x \to\infty} b(x) x^{1 - \alpha} > - \infty.\vspace*{-1pt}
\end{equation}
Then there exist positive constants $\gamma_0$ and $C_k(T)$, $k \ge
3$, such that\vspace*{-1pt}
%
%
\begin{equation}\label{e:expBound}
p_t(x, y)
\leq C_3(T)
\biggl( 1 + \frac{1}{ t^{3/2} } \biggr)
\exp\biggl( - \gamma_0 \frac{ ( y - x )^{ 2 (1 - \alpha) } }{ 2 C t
} \biggr)\vspace*{-1pt}
\end{equation}
and\vspace*{-1pt}
%
%
\begin{equation}\label{e:expDerivBound}
p^{(k)}_t(x, y)
\leq C_k(T)
\biggl( 1 + \frac{1}{ t^{ (2k + 3) / 2 } } \biggr)
\exp\biggl( - \gamma_0 \frac{ ( y - x )^{ 2 (1 - \alpha) } }{ 2 C t
} \biggr)\vspace*{-1pt}
\end{equation}
for every $y > x + 1$,
with $C = 2^{3 - 2\alpha} + 2|\gamma|_0^2 (1-\alpha)^2$.
The $C_k(T)$ also depend on $\alpha$ and on the coefficients $a, b$
and $\gamma$.\vspace*{-2pt}
\end{prop}
\begin{rem}
In the case of constant coefficients and $a = b = 0$, the bound (\ref
{e:expBound}) can be
compared to the density of the CEV process as provided, for example, in
\cite{Ford}, Theorem 1.6
(see also the references therein).
The comparison shows that our estimate is in the good range on the log-scale.\vspace*{-2pt}
\end{rem}
\begin{pf*}{Proof of Proposition \ref{p:asInfty}}
In the spirit of Lamperti's change-of-scale argument (cf. \cite{KS},
page 294), let $\varphi\in C^{2}( (0, \infty) )$ be
defined by
\[
\varphi(x)
= \int_0^x \frac{1}{ | \gamma|_0 y^{\alpha} } \,dy
= \frac{1}{ | \gamma|_0 (1 - \alpha) } x^{1 - \alpha},
\]
so that $\varphi'(x) \!=\! \frac{1}
{| \gamma|_0 x^{ \alpha} }$.
Let moreover $\theta\in C^{\infty}_b(\R)$ be such that \mbox{$1_{ [2,
\infty) } \!\le\!
\theta\!\le\!1_{ [1, \infty) }$} and $\theta' \le1$.
We set
\[
\rho(x)
= \cases{
\theta(x) \varphi(x), &\quad $x > 0$,
\cr
0, &\quad $x \le0$,}
\]
so that $\rho$ is of class $C^2 (\R)$.
We define the auxiliary process $Y_t = X_t - x$, $t \ge0$, which is
such that
$\Prob(Y_t \ge-x) = 1$, $t \ge0$.
An application of It\^o's formula yields
\[
\rho( Y_t )
=
\int_0^t f(Y_s) \,ds
+
M_t,
\]
where
\[
f(y)
= \rho'(y) \bigl( a(y) - b(y) y\bigr) + \tfrac{1}{2} \rho''(y) \gamma(y)^2
y^{2 \alpha}
\]
and
\[
M_t
= \int_0^t
\rho'(Y_s) \gamma(Y_s) Y_s^{\alpha} \,dW_s.
\]
The key point is the fact that $f$ is bounded from above on $(-x,
\infty)$ and,
on the other hand, $M$ is a martingale with bounded quadratic variation.
Indeed, $f$ is continuous, it is zero for $y \le1$ and for $y > 2$ one has
\[
f(y)
= \frac{ a(y) }{ |\gamma|_0 y^{ \alpha} }
- \frac{ b(y) }{ |\gamma|_0 } y^{ 1 - \alpha}
- \frac{\alpha}{2} \frac{ \gamma(y)^2 }{ |\gamma|_0 } y^{ \alpha-
1 },
\]
and hence, recalling that $a$ is bounded, $\varlimsup_{y \to\infty}
f(y) < \infty$
is ensured by
condition (\ref{e:bCondition1}).
Then we set $C_1 = \sup_{y \ge0} f(y)$.
For $M$, one has
%
%
\begin{eqnarray}\label{e:quadrVar}
\langle M \rangle_t
&=& \int_0^t \rho'(Y_s)^2 \gamma(Y_s)^2 Y_s^{2 \alpha} \,ds
\nonumber\\
&\leq&
\int_0^t \bigl( \varphi(2) + \varphi'(Y_s) \bigr)^2 \gamma(Y_s)^2
Y_s^{2 \alpha} \,ds
\\
&\leq&
2 \bigl( \varphi(2)^2 + 1 \bigr) t,\nonumber
\end{eqnarray}
and hence we set $C_2 = 2 ( \varphi(2)^2 + 1 ) =
2 ( \frac{2^{2(1-\alpha)}}{ | \gamma|_0^2 ( 1 - \alpha)^2 }
+ 1 )$.
Now, since $\rho$ is strictly increasing,\vspace*{2pt} $\{Y_t > y \} = \{ \rho
(Y_t) > \rho(y) \}$ for any $y > 0$.
Moreover,\vspace*{2pt}
$
\{ \rho(Y_t) > \rho(y) \}
\subseteq\{ M_t + C_1 t > \rho(y) \}
\subseteq\{ 2 M_t^2 + 2 C_1^2 t^2 > \rho(y)^2 \}
$.
We set\vspace*{2pt} $I_t = [(t-1) \vee t/2, t]$ and $\tau= \inf\{ s \ge0 \dvtx Y_s
\geq3/2 \}$.\vspace*{2pt}
The quadratic variation of $M$ is strictly increasing after $\tau$, since
$( \rho( Y_{\tau\vee t} ) \gamma( Y_{\tau\vee t} ) Y_{\tau\vee
t}^{\alpha} )^2 > 0$:
hence, by Dubins and Schwarz's theorem (cf. Theorem 3.4.6 in \cite
{KS}) there exists a one-dimensional Brownian motion $(b_t; t \ge0)$
such that $M_{\tau\vee t} = b_{ \langle M \rangle_{\tau\vee t} }$.
Clearly one has $\{ Y_t > 2 \} \subseteq\{ \tau< t \}$, so that for
$y > 2$
\begin{eqnarray*}
\overline{P}_t(y) &=& \Prob(\exists s \in I_t \dvtx Y_s > y)
\le
\Prob(\exists s \in I_t \dvtx Y_s > y, \tau< s)
\\[-1pt]
&\le&
\Prob\bigl(\exists s \in I_t \dvtx2 M_s^2 + 2 C_1^2 s^2 > \rho(y)^2, \tau< s\bigr)
\\[-1pt]
&\le&
\Prob\Bigl( \sup_{\tau< s \le t} ( 2 M_s^2 + 2 C_1^2 s^2
) > \rho(y)^2 \Bigr)
\\[-1pt]
&\le&
\Prob\biggl( \sup_{0 < s \le t} b_{ \langle M \rangle_s }^2 + C_1^2
t^2 > \frac{1}{2} \rho(y)^2 \biggr)
\\[-1pt]
&\le&
\Prob\biggl( \sup_{s \le C_2 t} b_s^2 + C_1^2 t^2 > \frac{1}{2} \rho
(y)^2 \biggr).
\end{eqnarray*}
%
We now employ the scaling property for the Brownian motion $( b_{s}; s
\ge0 ) \sim
( \sqrt{a} b_{s/a}; s \ge0 )$, $a > 0$, and Fernique's Theorem
(cf. \cite{IW}, page 402).
The latter tells that there exists a positive constant $\gamma_0$
such that $\exp( \gamma_0 \sup_{s \le1} b_s^2 )$ is integrable, hence
\begin{eqnarray*}
\overline{P}_t(y)
&\leq&\Prob\biggl( \gamma_0 \sup_{s \le1} b_s^2
+ \gamma_0 \frac{ C_1^2 }{ C_2 } t > \frac{\gamma_0}{2 C_2 t} \rho
(y)^2 \biggr)
\\
&\leq&
\exp\biggl( - \gamma_0 \frac{ \rho(y)^2 }{ 2 C_2 t } \biggr)
\esp[ e^{ \gamma_0 { C_1^2 }/{ C_2 } t + \gamma_0
\sup_{s \le1} b_s^2 } ]
\\
&\leq&
C_0
\exp
\biggl(
- \gamma_0 \frac{ \rho(y)^2 }{ 2 C_2 t }
+ \gamma_0 \frac{ C_1^2 }{ C_2 } t
\biggr),
\end{eqnarray*}
where $C_0 = \esp[ \exp( \gamma_0 \sup_{s \le1} b_s^2 ) ]$ is a
universal constant. The estimates on the density of $X_t$ and its
derivatives now follow from Theorem \ref{t:localDens} [estimates~%
(\ref{e:densEstimate}) and (\ref{e:densDerivEstimate})] and Theorem
\ref{t:localDensSubLin}(b), using $X_t - x = Y_t$, the value of the
constant $C_2$ and taking, for example, $R = 1/6$.
\end{pf*}

\subsection{Asymptotics at $0$} \label{s:zero}

We have established conditions under which the solution of (\ref
{e:localCIR}) admits a smooth density $p_t$ on $(0, + \infty)$.
According to Proposition~\ref{p:localCIR}, the process remains almost
surely in $\R_+$: this trivially means that for any $t > 0$, $X_t$ has
an identically zero density on $(- \infty, 0)$, which can be extended
to $0$ when $\tau_0 = \infty$ a.s. We are now wondering what are
sufficient conditions for $p_t$ to converge to zero at the origin,
hence providing the existence of a continuous (eventually
differentiable, eventually $C^{\infty}$) density on the whole real
line.

What we have in mind is the application of Theorem \ref
{t:localDensSubLin} to the inversed process
$Y_t = \frac{1}{X_t}$ (considered on the event $\{ \tau_0 = \infty\}$).
An application of It\^o's formula yields
%
%
\begin{equation}\label{e:inversCIR}
d Y_t
= J_{\alpha} ( Y_t )\,dt - \gammaHat(Y_t) Y_t^{ 2 - \alpha} \,dW_t,
\end{equation}
where
\[
J_{\alpha} (Y_t ) = - \widehat{a} (Y_t) Y_t^2
+ \widehat{b} (Y_t) Y_t
+ \gammaHat(Y_t)^2 Y_t^{ 3 - 2 \alpha}
\]
with the notation $\widehat{f} (y) = f(1/y)$, $y > 0$, for $f = a, b,
\sigma$.
Equation (\ref{e:inversCIR}) has super-linear coefficients, in
particular condition (\ref{e:polGrowth})
of Theorem \ref{t:localDensSubLin1} holds with \mbox{$q = 2$}.
Willing to apply Theorem \ref{t:localDensSubLin}(c), we first need
some preliminary results on the
moments of $Y$.
The proof of the next statement is based on the techniques employed in
\cite{BD2}, proof of Lemma
2.1, that we adapt to our framework.
\begin{lem} \label{l:invCIRmoments}
\begin{enumerate}[(1)]
\item[(1)] If $\alpha> 1/2$, assume \textup{(s1)$'$}.
Then for any initial condition $x > 0$, for any $t > 0$ and $p > 0$
%
%
\begin{equation}\label{e:invCIRMoments}
\esp
\biggl[
\sup_{s \le t}
\frac{1}{X_s^p}
\biggr]
\leq C,
\end{equation}
for some positive constant $C$ depending on $x, p, \alpha, t$ and on
the coefficients of~(\ref{e:localCIR}).
\item[(2)] If $\alpha= 1/2$, then assume \textup{(s1)} and
\textup{(s2)} and let
%
%
\begin{equation}\label{e:lStar}
l^* = \varliminf_{x \to0} \frac{ 2 a(x)} {\gamma(x)^2} > 1.
\end{equation}
Then (\ref{e:invCIRMoments}) holds for $p > 0$ such that
%
%
\begin{equation}\label{e:strongCIRcondition}
p + 1 < l^*.
\end{equation}
\end{enumerate}
\end{lem}
\begin{pf}
Let $\tau_n$ be the stopping time defined by $\tau_n = \inf\{ t \ge
0\dvtx X_t \le1/n \}$.
The application of It\^o's formula to $X_{t \wedge\tau_n}^p$, $p >
0$, yields
%
%
\begin{equation}\label{e:invCIRMoments1}
\esp \biggl[ \frac{1}{X_{t \wedge\tau_n}^p} \biggr] = \frac{1}{x^p}+\mathbb E \int_0^{t\wedge \tau_n}\varphi(X_s)\,ds,
\end{equation}
where
\[
\varphi(x) = p \frac{ b(x) }{ x^p } + \frac{ p }{ x^{p + 1} }
\underbrace{\biggl( \frac{ p + 1 }{ 2 } \gamma(x)^2 x^{ 2 \alpha- 1 } -
a(x) \biggr)}_{ g(x) },\qquad x > 0.
\]
It is easy to see that, if
%
%
\begin{equation}\label{e:gCondition}
\varlimsup_{x \to0} g(x)
< 0,
\end{equation}
there exists a positive constant $C$ such that $\varphi(x) < p \frac{
| b |_0 }{ x^p } + C$,
for every $x > 0$.
If (\ref{e:gCondition}) holds, from (\ref{e:invCIRMoments1}) we get
\[
\esp \biggl[ \sup_{s \le t} \frac{1}{X_{s \wedge\tau_n}^p} \biggr] \leq
\frac{1}{x^p} + C t + p | b |_0 \int_0^t \esp \biggl[ \sup_{u \leq s}
\frac{ 1 }{ X_{u \wedge\tau_n}^p } \biggr] \,ds,
\]
and hence, by Gronwall's lemma,
%
%
\begin{equation}\label{e:almostDoneGronw}
\esp \biggl[ \sup_{s \le t}\frac{1}{X_{s \wedge\tau_n}^p} \biggr] \leq \biggl(
\frac{1}{x^p} + C t \biggr) e^{p | b |_0 t}.
\end{equation}
We verify (\ref{e:gCondition}), distinguishing the two cases.

\textit{Case $\alpha> 1/2$.}
We simply observe that $\lim_{x \to0} g(x) = - a(0) < 0$.
Estimate~(\ref{e:invCIRMoments}) then follows by taking the limit $n
\rightarrow\infty$ in~(\ref{e:almostDoneGronw})
and using~Propo\-sition~\ref{p:localCIR} under assumption (s1)$'$.

\textit{Case $\alpha= 1/2$.}
We have $g(x) = \frac{p + 1}{ 2 } \gamma^2(x) - a(x)$.
If $p$ satisfies (\ref{e:strongCIRcondition}), then~(\ref{e:lStar})
ensures that $\varlimsup_{x \to0} g(x) < 0$.
We conclude again taking the limit $n \to\infty$ and using
Proposition \ref{p:localCIR} under assumptions (s1)
and (s2).
\end{pf}

We are now provided with the tools to prove the following:
\begin{prop} \label{p:asZero}
\begin{enumerate}[(1)]
\item[(1)] If $\alpha> 1/2$, assume \textup{(s1)$'$}. Then for every $t
> 0$, every $p \ge0$ and every $k > 0$ the density $p_t$ of $X_t$ on
$(0, +\infty)$ is such that
%
%
\begin{eqnarray}\label{e:contDens}
\lim_{y \rightarrow0^+} y^{- p} | p_t (y) | &=& 0,
\nonumber\\[-8pt]\\[-8pt]
\lim_{y \rightarrow0^+} y^{- p} \bigl| p^{(k)}_t (y) \bigr| &=& 0. \nonumber
\end{eqnarray}
\item[(2)] If $\alpha= 1/2$, then assume \textup{(s1)} and
\textup{(s2)} and define $l^*$ as in Lemma
\ref{l:invCIRmoments}.~If
%
%
\begin{equation}\label{e:conditionContDens}
l^* > 3 + q'_3 (2)
\end{equation}
[where $q'_{\cdot}(\cdot)$ has been defined in Theorem \ref{t:localDensSubLin}], then
%
%
\begin{equation}\label{e:attachment}
\lim_{y \rightarrow0^+} y^{- p} p_t (y) = 0
\end{equation}
for every $0 \leq p < l^* - ( 3 + q'_3 (2) )$.
Moreover, if
%
%
\begin{equation}\label{e:conditionContDerivDens}
l^* > 2 k + 3 + q'_{2 k + 3} (2),
\end{equation}
then
%
%
\begin{equation}\label{e:contDerivDens}
\lim_{y \rightarrow0^+} y^{- p} \bigl| p^{(k)}_t (y) \bigr| = 0
\end{equation}
for every $0 \le p < l^* - ( 2 k + 3 + q'_{2 k + 3} (2) )$.
\end{enumerate}
\end{prop}
\begin{pf}
We apply Theorem \ref{t:localDensSubLin1} to $Y = 1/X$. For simplicity
of notation, we write $p$ for $p_t$ and $p_Y$ for the density of $Y_t$.
As $Y$ satisfies (\ref{e:inversCIR}), from Theorem~\ref
{t:localDensSubLin1}(b) it follows that the bound (\ref{e:LambdaBound})
on $\Lambda_k$ holds with $q_k'(q) = q'_k(2)$. Hence, from Theorem
\ref{t:localDensSubLin1}(c) it follows that
%
%
\begin{equation}\label{e:limitInvDens}
\lim_{y' \rightarrow+ \infty}
p_Y (y') |y'|^{p'} = 0
\end{equation}
and, for a given $ k \geq0$,
%
%
\begin{equation}\label{e:limitInvDerviDens}
\lim_{y' \rightarrow+ \infty} p^{(k)}_{Y} (y') |y'|^{p'} = 0,
\end{equation}
if $\sup_{s \le t} Y_s$ has a finite moment of order $r > p' + q'_{2 k
+ 3} ( 2 )$.

Now, it is easy to see that
%
%
\begin{equation}\label{e:invDens}
p(y) = \frac{1}{y^2} p_{Y} \biggl( \frac{1}{y} \biggr),
\end{equation}
and hence, after some rather straightforward computations,
%
%
\begin{equation}\label{e:derivAndInvDeriv}
\bigl| p^{(k)} (y) \bigr| \leq C_k \biggl( \frac{1}{y} \biggr)^{2 (k +
1) } \sum_{j = 0}^k \sum_{\nu= 1}^j \frac{ d^{\nu} }{ dy^{\nu} } p_{Y}
\biggl( \frac{1}{y} \biggr),\qquad 0 < y < 1.
\end{equation}
Once again, we distinguish the two cases.

\textit{Case $\alpha> 1/2$.} If $1/2 < \alpha< 1$, by Lemma \ref
{l:invCIRmoments}, (\ref{e:limitInvDens}) and
(\ref{e:limitInvDerviDens}) hold for any $p' > 0$. Then
(\ref{e:contDens}) easily follows from (\ref{e:invDens}) and
(\ref{e:derivAndInvDeriv}).

\textit{Case $\alpha= 1/2$.} By Lemma \ref{l:invCIRmoments}, (\ref
{e:conditionContDens}) is the condition for $\sup_{s \le t} Y_s$
to have finite moment of order strictly greater than $2 + q'_3 ( 2 ) +
p$, with $p < l^* - ( 3 + q'_3 (2) )$.
By (\ref{e:invDens}), in this case (\ref{e:limitInvDens}) holds true
with $k = 0$ and $p' = 2 + p$, and hence
(\ref{e:attachment}) holds.
Similarly, by (\ref{e:derivAndInvDeriv}) estimate (\ref
{e:contDerivDens}) holds if
(\ref{e:limitInvDens}) holds with $p' = p + 2 (k + 1)$.
The latter condition is achieved if $\sup_{s \le t} Y_s$ has finite
moment of order strictly greater than
$2 (k + 1) + q'_{2 (k + 3)} ( 2 ) + p$, which is in turn ensured by
(\ref{e:conditionContDerivDens}).
\end{pf}
\begin{rem}
Proposition \ref{p:asInfty} states that $p_t$ decays exponentially at
infinity for any value of $\alpha$,
as far as condition (\ref{e:bCondition1}) holds true.
When $\alpha> 1/2$, Proposition~\ref{p:asZero} states that $p_t$ and
all its derivatives tend to zero at the origin, while the price to pay
for the same conclusion to hold is higher when $\alpha=
1/2$ [cf. conditions (\ref{e:conditionContDens}) and (\ref
{e:conditionContDerivDens}), which become rapidly strong
for growing values of $k$].
With regard to this behavior at zero, we recall that Proposi\-tion~\ref
{p:asZero} only provides sufficient
conditions for estimates (\ref{e:attachment}) and (\ref
{e:contDerivDens}) to hold.
We do not give any conclusion on the behavior of the density at zero
when condition (\ref{e:conditionContDens}) [or (\ref
{e:conditionContDerivDens}) for the derivatives] fail to hold.
\end{rem}

\begin{appendix}
\section*{Appendix}

We collect here the proofs of some of the more technical results.

\subsection{\texorpdfstring{Proofs of Lemmas \protect\ref{l:boundDeriv} and
\protect\ref{l:covMatrix}}{Proofs of Lemmas 2.1 and 2.2}} \label{a:app1}

\mbox{}

\begin{pf*}{Proof of Lemma \ref{l:boundDeriv}}
We refer to the notation introduced in the proof of~\cite{Nual},
Theorem 2.2.2,
allowing us to write the equation satisfied by the $k$th
Malliavin derivative in a compact form.
This is stated as follows: for any subset $K = \{h_1,\ldots,h_{\eta}\}
$ of $\{1,\ldots,k\}$,
one sets $j(K) = j_{h_1},\ldots, j_{h_{\eta}}$ and\vadjust{\eject} $r(K) =
r_{h_1},\ldots, r_{h_{\eta}}$.
Then, one defines
%
\begin{eqnarray*}
\alpha^i_{l,j_1,\ldots,j_k}(s,r_1,\ldots, r_k)
:\!&=& D^{j_1,\ldots, j_k}_{r_1,\ldots, r_k} A^i_l( X_s )
\\
&=& \sum\partial_{k_1} \cdots\partial_{k_{\nu}} A^i_l(X_s)
\\
&&\hspace*{12.5pt}{} \times D^{j(I_1)}_{r(I_1)} X^{k_1}_s \cdots D^{j(I_{\nu
})}_{r(I_{\nu})} X^{k_{\nu}}_s
\end{eqnarray*}
and
\begin{eqnarray*}
\beta^i_{j_1,\ldots,j_k}(s,r_1,\ldots, r_k)
:\!&=& D^{j_1,\ldots, j_k}_{r_1,\ldots, r_k} B^i( X_s )
\\
&=& \sum\partial_{k_1} \cdots\partial_{k_{\nu}} B^i(X_s)
\\
&&\hspace*{12.5pt}{} \times D^{j(I_1)}_{r(I_1)} X^{k_1}_s \cdots D^{j(I_{\nu
})}_{r(I_{\nu})} X^{k_{\nu}}_s,
\end{eqnarray*}
where in both cases the sum is extended to the set of all partitions of
$\{1,\ldots,k\} =
I_1 \cup\cdots\cup I_{\nu}$.
Finally, one sets $\alpha^i_j(s) = A^i_j(X_s)$.
Making use\vspace*{1pt} of this notation, it is shown that the equation satisfied by
the $k$th derivative
reads as
%
%
\setcounter{equation}{0}
\begin{eqnarray} \label{e:linearEq}\quad
\DX&=& \sum_{\varepsilon=1}^k \alpha^i_{j_{\varepsilon},j_1,\ldots
,j_{\varepsilon-1},j_{\varepsilon+1},\ldots,j_k}
(r_{\varepsilon},r_1,\ldots,r_{\varepsilon-1},r_{\varepsilon+1},\ldots,r_k)
\nonumber\\
&&{} +\int_{r_1 \vee\cdots\vee r_k}^t \bigl( \beta^i_{j_1,\ldots
,j_k}(s,r_1,\ldots, r_k) \,ds\\
&&\hspace*{53.2pt}{}
+ \alpha^i_{l,j_1,\ldots,j_k}(s,r_1,\ldots, r_k) \,dW^l_s
\bigr),\nonumber
\end{eqnarray}
if $t \geq r_1 \vee\cdots\vee r_k$, and $\DX= 0$ otherwise. We prove
(\ref{e:boundDeriv}) by induction. The estimate is true for $k = 1$,
with $\gamma_{1, p} = 2 C_{1,p}$: this simply follows with an
application of Burkholder's inequality and Gronwall's lemma to (\ref
{e:linearEq}) taken for $k = 1$. Let us suppose that
(\ref{e:boundDeriv}) is true up to $k-1$. As done for $k = 1$, we apply
Burkholder's inequality to (\ref {e:linearEq}), and, setting $r = r_1
\vee\cdots\vee r_k$, we get
%
%
\begin{eqnarray} \label{e:e1}
&&\esp[ | \DX|^p ]\nonumber\\
&&\qquad\leq C_{k,p}
\Biggl\{ \sum_{\varepsilon=1}^k \esp[ |
\alpha^i_{j_{\varepsilon},j_1,\ldots,j_{\varepsilon-1},j_{\varepsilon
+1},\ldots,j_k}
(r_{\varepsilon},r_1,\ldots,r_{\varepsilon-1},r_{\varepsilon+1},\ldots
,r_k) |^p ]
\nonumber\\
&&\hspace*{24.1pt}\qquad\quad{}+ (t - r)^{{p}/{2}-1}\nonumber\\
&&\qquad\quad\hspace*{34.1pt}{}\times \mathop{\sum_{I_1 \cup\cdots\cup
I_{\nu }}}_{\mathrm{card}(I) \leq k-1} \int_r^t \esp \Biggl[ \Biggl(
(t-r)^{{1/2}} | \partial_{k_1} \cdots\partial_{k_{\nu}}
B^i(X_s) |\nonumber\\
&&\hspace*{24.1pt}\qquad\quad\hspace*{121.1pt}{}
+ \sum_{l=1}^d | \partial_{k_1} \cdots\partial_{k_{\nu}}
A^i_l(X_s) |
\Biggr)^p
\\
&&\hspace*{24.1pt}\qquad\quad\hspace*{112.7pt}{}\times\bigl|
D^{j(I_1)}_{r(I_1)} X^{k_1}_s \cdots D^{j(I_{\nu})}_{r(I_{\nu })}
X^{k_{\nu}}_s \bigr|^p \Biggr] \,ds
\nonumber\\
&&\hspace*{24.1pt}\qquad\quad{}
+ (t - r)^{{p}/{2}-1} \int_r^t
\esp\Biggl[ \Biggl(
(t-r)^{{1/2}} | \partial_k B^i(X_s) |\nonumber\\
&&\hspace*{24.1pt}\qquad\quad\hspace*{122.6pt} + \sum_{l=1}^d |
\partial_k A^i_l(X_s) | \Biggr)^p\nonumber\\
&&\hspace*{86.3pt}\qquad\quad\hspace*{112.7pt}{}\times
| D^{j_1,\ldots,j_k}_{r_1,\ldots,r_k}
X_s^k |^p \Biggr] \,ds \Biggr\},\nonumber
\end{eqnarray}
where, in the last line, we have isolated the term depending on
$D^{j_1,\ldots,j_k}_{r_1,\ldots,r_k} X$.

To estimate the second term in (\ref{e:e1}) we notice that, for any
partition $I_1 \cup\cdots\cup I_{\nu}$
of $\{1,\ldots, k \}$ such that $\mathrm{card}(I) \leq k-1$, using
(\ref{e:boundDeriv}) up to order
$k-1$ we have
%
%
\begin{eqnarray} \label{e:lambda1}
&&\esp\Biggl[
\Biggl(
(t-r)^{{1/2}} | \partial_{k_1} \cdots\partial_{k_{\nu}}
B^i(X_s) |
\nonumber\\
&&\qquad\quad\hspace*{3.5pt}{}+\sum_{l=1}^d | \partial_{k_1}
\cdots\partial_{k_{\nu}} A^i_l(X_s) | \Biggr)^p\bigl |
D^{j(I_1)}_{r(I_1)} X^{k_1}_s \cdots D^{j(I_{\nu})}_{r(I_{\nu})}
X^{k_{\nu}}_s \bigr|^p \Biggr]
\\
&&\qquad\leq
C \{
\nA^{k}_{k - 2}
( t^{1/2} \nB_k + \nA_k)^{\chi_k} \}^p
e_p(t)^{ \lambda^{(1)}_{k, p} },\nonumber
\end{eqnarray}
where we have defined
\[
\lambda^{(1)}_{k, p} := \mathop{\sup_{I_1 \cup\cdots\cup I_{\nu
} = \{1,\ldots,k\}}}_{\mathrm{card}(I) \leq k-1} \bigl\{
\gamma_{\mathrm{card}(I_1), p} + \cdots+ \gamma_{\mathrm{card}(I_{\nu
}), p}
\bigr\}
\]
and
\[
\chi_k
= 1 + \sum_{l=1}^{\nu} \bigl( \mathrm{card}(I_l) + 1\bigr)^2.
\]
It is easy to see that
\[
\chi_k \leq(k + 1)^2,
\]
since
\begin{eqnarray*}
\sum_{l=1}^{\nu} \mathrm{card}(I_l)^2
=
(k-1)^2 \sum_{l=1}^{\nu} \frac{ \mathrm{card}(I_l)^2 }{ (k -1)^2 }
\leq
(k-1)^2 \sum_{l=1}^{\nu} \frac{ \mathrm{card}(I_l) }{ (k -1) }
=
(k-1)k,
\end{eqnarray*}
so that
\begin{eqnarray*}
\chi_k &=& 1 + \sum_{l=1}^{\nu} \mathrm{card}(I_l)^2 + 2
\sum_{l=1}^{\nu} \mathrm{card}(I_l) + \nu
\\
&\leq& 1 + ( k - 1 ) k + 2 k + \nu
\\
&\leq& 1 + k^2 - k + 3 k = ( k +1 )^2.
\end{eqnarray*}
To estimate the first term in (\ref{e:e1}), notice that we have as well
%
%
\begin{eqnarray} \label{e:lambda2}
&&\esp[ |
\alpha^i_{j_{\varepsilon},j_1,\ldots,j_{\varepsilon-1},j_{\varepsilon
+1},\ldots,j_k}
(r_{\varepsilon},r_1,\ldots,r_{\varepsilon-1},r_{\varepsilon+1},\ldots,r_k)
|^p ]
\nonumber\\[-8pt]\\[-8pt]
&&\qquad\leq C \{ \nA^{k}_{k-1} ( r_{\varepsilon}^{{1/2}} \nB_{k-1} +
\nA_{k-1} )^{k^2} \}^p e_p(t)^{ \lambda^{(2)}_{k, p} },\nonumber
\end{eqnarray}
with
\[
\lambda^{(2)}_{k, p} := \sup_{I_1 \cup\cdots\cup I_{\nu} = \{
1,\ldots,k-1\}} \bigl\{
\gamma_{\mathrm{card}(I_1), p} + \cdots+ \gamma_{\mathrm{card}(I_{\nu
}), p}
\bigr\}.
\]
We remark that $\lambda^{(1)}_k$ and $\lambda^{(2)}_k$ are defined by
means of the $\gamma$'s
up to order $k-1$.

Collecting (\ref{e:e1}), (\ref{e:lambda1}) and (\ref{e:lambda2}), we get
\begin{eqnarray*}
&&
\esp[ | \DX|^p ] \\
&&\qquad\leq C_{k,p} \Biggl\{ \nA^{kp}_{k-1} ( t^{1/2} \nB_{k-1} +
\nA_{k-1} )^{k^2 p} e_p(t)^{ \lambda^{(2)}_{k,p} }
\\
&&\qquad\quad\hspace*{26pt}{} + (t-r)^{{p/2}} \nA^{k}_{k-2} ( t^{1/2} \nB_k + \nA_k )^{( k
+1 )^2 p} e_p(t)^{ \lambda^{(1)}_{k, p} }
\\
&&\qquad\quad\hspace*{26pt}{} + (t-r)^{{p/2}-1} ( t^{1/2} \nB_1 + \nA_1 )^p \sum_{k=1}^m
\int_r^t \esp| D^{j_1,\ldots, j_k}_{r_1,\ldots,r_k} X^i_s |^p \,ds \Biggr\}
\\
&&\qquad\leq C_{k,p} \nA^{kp}_{k-1} ( t^{1/2} \nB_k + \nA_k )^{( k +1 )^2 p}
e_p(t)^{ \lambda^{(1)}_{k, p} \vee\lambda^{(2)}_{k, p} }
\\
&&\qquad\quad{} \times\bigl( 1 + C_{k,p} t^{p/2} ( t^{1/2} \nB_1 + \nA_1 )^p e_p(t)^{
C_{k,p} } \bigr)
\\
&&\qquad\leq C_{k,p} \nA^{kp}_{k-1} ( t^{1/2} \nB_k + \nA_k )^{( k +1 )^2 p}
e_p(t)^{ \lambda^{(1)}_{k, p} \vee\lambda^{(2)}_{k, p} + 2 C_{k,p} },
\end{eqnarray*}
where we have applied Gronwall's lemma to get the second inequality.
The constant $C_{k,p}$ may vary from line to line, but never depends on
$t$ nor on the
bounds on $B$ and $A$.
We recursively define $\gamma_{k, p}$ by setting $\gamma_{k, p} :=
\lambda^{(1)}_{k, p} \vee
\lambda^{(2)}_{k, p} + 2 C_{k, p}$,
and we finally obtain (\ref{e:boundDeriv}).
\end{pf*}
\begin{pf*}{Proof of Lemma \ref{l:covMatrix}}
\textit{Step} 1. We first use the decomposition $D_s X_t = Y_t Z_s
A(X_s)$ (see, e.g., \cite{Nual}) and write
%
%
\begin{eqnarray} \label{E:sigma1}
\s&=& Y_t \int_0^t Z_s A(X_s) A(X_s)^* Z_s^* \,ds\, Y_t^*
\nonumber\\[-8pt]\\[-8pt]
&=& Y_t U_t Y_t^*,
\nonumber
\end{eqnarray}
where we have set $U_t = \int_0^t Z_s A(X_s) A(X_s)^* Z_s^* \,ds$.
Notice that $U_t$ is a positive operator,
and that for any $\xi\in\R^m$ we have
\[
\langle  \xi, U_t \xi\rangle  = \int_0^t \langle A(X_s)^* Z_s^* \xi, A(X_s)^* Z_s^* \xi
\rangle  \,ds =
\sum_{j=1}^d \int_0^t \langle  Z_s A_j(X_s), \xi\rangle ^2.
\]
From identity (\ref{E:sigma1}) it follows that $\det\s=
(\det Y_t)^2 \det U_t = (\det Z_t)^{-2} \det U_t$.
Hence, applying H\"older's inequality,
%
%
\begin{eqnarray}\label{E:sigma2}
\mathbb{E}[|{\det\s}|^{-p}] &\leq& ( \esp[ |{\det Z_t}|^{4p} ] \esp[
(\det U_t)^{-2p} ] )^{1/2}
\nonumber\\[-8pt]\\[-8pt]
&\leq& C_{p,m} ( e^Z_{4p}(t)^m \esp[ (\det U_t)^{-2p} ] )^{1/2},
\nonumber
\end{eqnarray}
where in the last step we have used bound (\ref{E:ZBound}) on the
entries of $Z_t$.

\textit{Step} 2.
Let $\lambda_t = \inf_{| \xi|=1} \langle \xi, U_t \xi\rangle $ be the smallest
eigenvalue of $U_t$, so
that $\mathbb{E}[ (\det U_t)^{ -2p }] \leq\mathbb
{E}[ \lambda_t^{-2mp}]$.
We evaluate $\Prob(\lambda_t \leq\varepsilon)$.

For any $\xi$ such that $| \xi| = 1$, using the elementary inequality
$(a+b)^2 \geq a^2/2 - b^2$
we get
\begin{eqnarray*}
\sum_{j=1}^d \langle  Z_s A_j(X_s), \xi\rangle ^2 &\geq&
\frac{1}{2} \sum_{j=1}^d \langle  A_j(x), \xi\rangle ^2
-
\sum_{j=1}^d \langle  Z_s A_j(X_s) - A_j(x), \xi\rangle ^2
\\
&\geq&\frac{1}{2} c_* - \sum_{j=1}^d | Z_s A_j(X_s) - A_j(x) |^2,
\end{eqnarray*}
where in the last step we have used the ellipticity assumption (E).
For any $\varepsilon> 0$ and $a > 0$ such that $a\varepsilon< t$, the
previous inequality gives
\[
\Prob(\lambda_t \leq\varepsilon) \leq \Prob \Biggl( \frac{1}{2} a c_*
\varepsilon - \sup_{s\leq a \varepsilon} \Biggl\{ a \varepsilon
\sum_{j=1}^d | Z_s A_j(X_s) - A_j(x) |^2 \Biggr\} \leq\varepsilon
\Biggr),
\]
and thus, if we take $a = 4/c_*$ in order to have $a c_*/ 2 = 2$ and
apply Markov's inequality,
we obtain
%
%
\begin{eqnarray}\label{E:lambda2}
\Prob(\lambda_t \leq\varepsilon) &\leq& \Prob \Biggl( \sup_{s\leq a
\varepsilon} \Biggl\{ \sum_{j=1}^d | Z_s A_j(X_s) - A_j(x) |^2 \Biggr\}
\geq\frac{c_*}{4} \Biggr)
\nonumber\\[-8pt]\\[-8pt]
&\leq& d^{q-1} \frac{4^q}{c_*^q} \sum_{j=1}^d
\esp\Bigl[ {\sup_{s\leq a \varepsilon}} |Z_s A_j(X_s) - A_j(x)|^{2q}
\Bigr],\nonumber
\end{eqnarray}
where the last holds for all $q > 1$.
Now, to estimate the last term we claim that, for all $j = 1,\ldots,
d$,\vspace*{-1pt}
%
%
\begin{equation}\label{E:ZAIncrem}
\esp \Bigl[ {\sup_{s\leq t}} |Z_s A_j(X_s) - A_j(x)|^{2q} \Bigr] \leq C
t^q ( t^{1/2}|B|_0|A|_2^3 + |A|_1^2 )^{2q} e^Z_{2q}
(t)^C,\hspace*{-28pt}\vspace*{-1pt}
\end{equation}
for a constant $C$ depending on $q, m, d$ but not on the bounds on $B$
and $A$.
From (\ref{E:lambda2}) and this last estimate, it follows that\vspace*{-1pt}
\[
\Prob(\lambda_t \leq\varepsilon) \leq
C_{q,m,d}
\frac{\varepsilon^q}{c_*^{2q}}
( t^{1/2}|B|_0|A|_2^3 + |A|_1^2 )^{2q}
e^Z_{2q} (t)^{ C_{q,m,d} },\vspace*{-1pt}
\]
for any $\varepsilon$ such that $4 \varepsilon/c_* < 1 \wedge t$.

\textit{Step} 3. We finally estimate $\mathbb{E}[\lambda_t^{-2mp}]$.
We write\vspace*{-1pt}
\begin{eqnarray*}
\mathbb{E}[\lambda_t^{-2mp}] &=&
\mathbb{E}\bigl[\lambda_t^{-2mp} 1_{\{\lambda_t > 1\}}\bigr] +
\sum_{k=1}^{\infty} \mathbb{E}\bigl[\lambda_t^{-2mp} 1_{\{ 1/(k+1) <
\lambda_t \leq1/k \}}\bigr] \\[-5pt]
&\leq&1 + \sum_{k=1}^{\infty} (k+1)^{2mp} \Prob\bigl( 1/(k+1) <
\lambda_t \leq1/k\bigr),\vspace*{-1pt}
\end{eqnarray*}
and separate the contribution of the sum over $k> \frac{4}{t c_*}$ to
obtain
\begin{eqnarray*}
\mathbb{E}[\lambda_t^{-2mp}]
&\leq&
1 + \sum_{1\leq k\leq{4}/({t c_*})} (k+1)^{2mp} \Prob\bigl( 1/(k+1) <
\lambda_t \leq1/k\bigr)
\\[-2pt]
&&{}
+ \sum_{k > {4}/({t c_*})} (k+1)^{2mp} \Prob(\lambda_t \leq1/k)
\\[-2pt]
&\leq& 1 + \biggl( \frac{4}{t c_*} + 1 \biggr)^{2mp}
\Prob( \lambda_t \leq1)
\\[-2pt]
&&{}
+ C_{q,m,d}
e^Z_{2q}(t)^{ C_{q,m,d} }
\frac{1}{c_*^{2q}}
( t^{1/2}|B|_0|A|_2^3 + |A|_1^2 )^{2q}
\\[-2pt]
&&\hspace*{10.3pt}{} \times
\sum_{k > {4}/({t c_*})} (k+1)^{2mp} \frac{1}{k^q}.
\end{eqnarray*}
We finally take $q = 2mp + 2$ in order to get convergent series.
This last estimate, together with (\ref{E:sigma2}), gives the desired
result.\vspace*{-3pt}\noqed
\end{pf*}

\begin{pf*}{Proof of (\ref{E:ZAIncrem})}
We apply It\^o's formula to the product $Z_t A_j(X_t)$~and~get\vspace*{-1pt}
%
%
\begin{eqnarray}\label{e:ZA}\qquad
d( Z_t A_j(X_t)) &=&
Z_t
\Biggl\{
(\dA_j B - \partial B A) + \sum_{l=1}^d \dA_l (\dA_l A_j - \dA_j A_l)
\Biggr\} (X_t) \,dt
\nonumber\\[-2pt]
&&{} + Z_t
\biggl(
\frac{1}{2} \partial_{k_1} \partial_{k_2} A_j A_l^{k_1} A_l^{k_2}
\biggr) (X_t) \,dt
\\[-2pt]
&&{} + Z_t \sum_{l=1}^d (\dA_j A_l - \dA_l A_j) (X_t) \,dW^l(t).\nonumber
\end{eqnarray}
Hence, by Burkholder's inequality,
\begin{eqnarray*}
&&\sup_{i=1,\ldots,m}
\mathbb{E}
\Bigl[
\sup_{s\leq t} \bigl| \bigl(Z_s A_j(X_s) - A_j(x)\bigr)^i \bigr|^{2q}
\Bigr]
\\
&&\qquad\leq C \{ t^{2q-1} (|B|_0|A|_1 + |B|_1|A|_0 + |A|_1^2|A|_0
+ |A|_2|A|_0^2)^{2q}
\\
&&\qquad\quad\hspace*{156.2pt}{}+ t^{q-1} ( |A|_1|A|_0 )^{2q} \}\\
&&\qquad\quad\hspace*{0pt}{}\times
\int_0^t
\esp\Bigl[ \sup_{h,k = 1,\ldots,m} |(Z_u)_{h,k}|^{2q} \Bigr] \,du
\\
&&\qquad\leq C t^q
(t^{1/2}|B|_0|A|_2^3 + |A|_1^2 )^{2q}
e^Z_{2q}(t)^C,
\end{eqnarray*}
where the constant $C$ depends on $q,m$ and $d$, but not on $t$
and on the bounds on $B$ and $A$ and their derivatives.
In the last step, we have once again used bound (\ref{E:ZBound}) on
the entries of $Z$.
\end{pf*}

\subsection{\texorpdfstring{Proof of Proposition \protect\ref{p:localCIR}}{Proof of Proposition 3.1}}
\label{a:app2}

We first collect the basic facts we need to give the proof of Proposition
\ref{p:localCIR}.
We will start by proving existence and uniqueness of strong solutions
for the following equation:
%
%
\begin{eqnarray}\label{e:localCIRmodule}
X_t &=& x
+ \int_0^t \bigl( a(X_s) - b(X_s) X_s \bigr) \,ds\nonumber\\[-8pt]\\[-8pt]
&&{} + \int_0^t \gamma(X_s) | X_s |^{\alpha} \,dW_s,\qquad t \ge0,
\alpha\in[1/2, 1),\nonumber
\end{eqnarray}
whose coefficients are defined on the whole real line [$a, b$ and
$\gamma$ are the
functions appearing in (\ref{e:localCIR})].
Once we have established that the unique strong solution of (\ref
{e:localCIRmodule}) is a.s. positive,
then (\ref{e:localCIRmodule}) will coincide with the original equation
(\ref{e:localCIR}).

The proof of Proposition \ref{p:localCIR} is split in the following
two short lemmas.
\setcounter{lem}{0}
\begin{lem} \label{l:localCIRmodule}
Assume condition \textup{(s0)} of Proposition \ref{p:localCIR}.
Then existence and uniqueness of strong solutions hold for (\ref
{e:localCIRmodule}).
Moreover, for any initial condition $x \ge0$ the solution is a.s.
positive, $\Prob(X_t \ge0; t \ge0)
= 1$.
\end{lem}
\begin{pf}
Existence of nonexplosive weak solutions for (\ref{e:localCIRmodule})
follows from continuity and
sub-linear growth of drift and diffusion coefficients.
The existence of weak solutions together with pathwise uniqueness imply
the existence of strong
solutions (cf. \cite{KS}, Proposition 5.3.20 and Corollary 5.3.23).
Pathwise uniqueness follows in its turn from a well-known theorem of
uniqueness of Yamada and
Watanabe (cf. \cite{KS}, Proposition 5.2.13).
Indeed, as $a, b, \gamma\in C^1_b$ the diffusion coefficient of (\ref
{e:localCIRmodule}) is
locally H\"older-continuous of exponent $\alpha\ge1/2$ and the drift
coefficient is locally
Lipschitz-continuous.
We apply the standard localization argument for locally Lipschitz
coefficients and the Yamada--Watanabe
theorem to establish that solutions are pathwise unique up to their
exit time from a compact ball, and hence
pathwise uniqueness holds for (\ref{e:localCIRmodule}).
\end{pf}

Lemma \ref{l:FellerCIR} deals with the second part of Proposition \ref
{p:localCIR}, that is, the
behavior at zero.
The proof is based on Feller's test for explosions of solutions of
one-dimensional SDEs (cf. \cite{KS},
Theorem 5.5.29).
Letting $\tau$ denote the exit time from $(0, \infty)$, that is, $\tau
= \inf\{ t \geq0 \dvtx X_t \notin
(0, \infty) \}$ with $\inf\varnothing= \infty$, we have to verify that
%
%
\begin{equation}\label{e:noExpl}
\lim_{x \to0} p_c(x) = -\infty
\end{equation}
with $p_c$ defined by
%
%
\begin{equation}\label{E:scaleFunct}
p_c(x) := \int_c^x
\exp{ \biggl( - 2 \int_c^y \frac{ a(z) - b(z) z }{ \gamma(z)^2
z^{2 \alpha} } \,dz \biggr) }\,
dy, \qquad x > 0,
\end{equation}
for a fixed $c > 0$.
Property (\ref{e:noExpl}) implies that $\Prob(\tau= \infty) = 1$, then
$\tau\equiv\tau_0$ with $\tau_0$ as defined in Proposition \ref
{p:localCIR},
because the solution of (\ref{e:localCIRmodule}) does not\vadjust{\goodbreak} explode
at $\infty$ (cf. Lemma \ref{l:localCIRmodule}).
The inner integral in (\ref{E:scaleFunct}) is well defined and finite
for any $y > 0$ because
$\gamma(z)^2 > 0$ for any $z > 0$ and $\gamma$ is continuous.
\setcounter{rem}{0}
\begin{rem}
The conclusion does not depend on the choice of \mbox{$c \in(0,\infty)$}.
\end{rem}
\begin{lem} \label{l:FellerCIR}
Assume \textup{(s0)}, and let $X = (X_t; t\ge0)$ denote the unique strong
solution of (\ref{e:localCIRmodule}) for initial condition $x > 0$.
Then the statements of Proposition~\ref{p:localCIR} on the stopping time
$\tau_0$ hold true.
\end{lem}
\begin{pf}
We prove (\ref{e:noExpl}), for $c = 1$.
We assume without restriction that $x < 1$ and distinguish the two cases.

\textit{Case $\alpha> 1/2$}.
We have $a(z) \geq a(0) - |a|_1 z$, $z > 0$.
Then
\begin{eqnarray*}
\frac{a(z) - b(z) z}{\gamma(z)^2 z^{2 \alpha}}
&\geq&\frac{a(0) - ( | a |_1 + b(z) ) z}{\gamma(z)^2 z^{2 \alpha}}\\
&\geq&\frac{a(0)}{| \gamma|_0^2 z^{2 \alpha}}
- \frac{ | a |_1 + | b |_0 }{\gamma(z)^2 z^{2 \alpha- 1}},
\end{eqnarray*}
$\frac{ 1 }{\gamma(z)^2 z^{2 \alpha- 1}}$ is integrable at zero by
(s1)$'$, and
then there exists a positive constant $K$ such that
\begin{eqnarray*}
- 2 \int_1^y \frac{a(z) - b(z) z}{\gamma(z)^2 z^{2 \alpha}}
&\geq&\frac{2 a(0)}{| \gamma|_0^2} \int_y^1 \frac{dz}{z^{2 \alpha}}
+ K\\
&=& \frac{2 a(0)}{(2 \alpha- 1) | \gamma|_0^2} \biggl( \frac{1}{y^{2
\alpha- 1}} - 1 \biggr) + K
\end{eqnarray*}
hence
\begin{eqnarray*}
p_1(x) &\le&
- C \int_x^1 \exp\biggl( \frac{2 a(0)}{(2 \alpha- 1) | \gamma|_0^2}
\frac{1}{y^{2 \alpha- 1}} \biggr) \,dy
\\
&=& - C \int_{ 1 }^{ {1/x} }
\frac{1}{t^2}
\exp\biggl( \frac{2 a(0)}{(2 \alpha- 1) | \gamma|_0^2} t^{2 \alpha
- 1} \biggr) \,dt
\mathop{\to}_{x \to0^+} - \infty.
\end{eqnarray*}

\textit{Case $\alpha= 1/2$}. By (s2),
\[
\frac{ 2 a(z) }{\gamma(z)^2 z } \geq\frac{1}{z},
\]
for $z < \overline{x}$.
Hence
\[
2 \frac{a(z) - b(z) z}{\gamma(z)^2 z}
\geq\frac{ 1 }{ z } - 2 | b |_0 \frac{1}{ \gamma(z)^2 },
\]
and thus, $\frac{1}{ \gamma^2 }$ being integrable at zero, for $x <
\overline{x}$
we have
\begin{eqnarray*}
p_1(x) &\leq&
- C \int_x^1 \exp\biggl( \int_y^1 \frac{ 1 }{ z } \,dz \biggr) \,dy
\\
&=& - C \int_x^1 \frac{1}{y} \,dy \mathop{\to}_{x \to0^+} - \infty.
\end{eqnarray*}
\upqed\end{pf}
\end{appendix}

\section*{Acknowledgments}
I am grateful to Vlad Bally of University Paris-Est Marne-la-Vall\'ee
for introducing me to this subject and for many stimulating discussions
and to Alexander Yu. Veretennikov of University of Leeds for useful
insights. I sincerely thank an anonymous referee for reading a previous
version of this paper and for providing helpful comments to improve the
presentation.

%

%
\printaddresses

\end{document}